\documentclass[twocolumn]{svjour3}
\usepackage[usenames,dvipsnames,svgnames,table]{xcolor}

\usepackage[utf8]{inputenc}
\usepackage[T1]{fontenc}
\usepackage{graphicx}
\usepackage{longtable}
\usepackage{wrapfig}
\usepackage{rotating}
\usepackage[normalem]{ulem}
\usepackage{amsmath}
\usepackage{amssymb}
\usepackage{capt-of}
\usepackage{hyperref}
\usepackage{algorithm}
\usepackage{algpseudocode}

\usepackage{pgf}
\usepackage{amsmath}

\newcommand{\rev}[2]{#2}

\DeclareMathOperator*{\argmin}{arg\,min}
\author{Mark Blyth \and Krasimira Tsaneva-Atanasova \and Lucia Marucci\(^*\) \and Ludovic Renson\(^*\) \newline \small * Co-last author
}
\institute{
Mark Blyth
  \at Department of Engineering Mathematics, University of Bristol, Bristol, UK.
  \and    Krasimira Tsaneva-Atanasova  \at Hub for Quantitative Modelling in Healthcare, University of Exeter, Exeter, UK.
  \at   The Alan Turing Institute, British Library, London, UK.
  \at    Data Science Institute, University of Exeter, Exeter, UK.
  \and Lucia Marucci  \at Department of Engineering Mathematics, University of Bristol, Bristol, UK.
  \at    BrisSynBio, University of Bristol, Bristol, U.K.
  \at    School of Cellular and Molecular Medicine, University of Bristol, Bristol, UK.
  \and Ludovic Renson \at  Department of Mechanical Engineering, Imperial College London, London, UK.
  \\\email{l.renson@imperial.ac.uk}
}
\journalname{Nonlinear Dynamics}
\date{}
\title{Numerical methods for control-based continuation of relaxation oscillations}

\begin{document}

\maketitle

\begin{abstract}
Control-based continuation (CBC) is an experimental method that can reveal stable and unstable dynamics of physical systems.
It extends the path-following principles of numerical continuation to experiments, and provides systematic dynamical analyses without the need for mathematical modelling.
CBC has seen considerable success in studying the bifurcation structure of mechanical systems.
Nevertheless, the method is not practical for studying relaxation oscillations.
Large numbers of Fourier modes are required to describe them, and \rev{experiments become impractically slow}{the length of the experiment significantly increases} when many Fourier modes are used\rev{}{, as the system must be run to convergence many times}.
Furthermore, relaxation oscillations often arise in autonomous systems, for which an appropriate phase constraint is required.
To overcome these challenges, we introduce an adaptive B-spline discretisation, that can produce a parsimonious description of responses that would otherwise require many Fourier modes.
We couple this to a novel phase constraint, that phase-locks control target and solution phase.
Results are demonstrated on simulations of a slow-fast synthetic gene network and an Oregonator model.
Our methods extend CBC to a much broader range of systems than have been studied so far, opening up \rev{}{a range of novel} experimental opportunities on slow-fast systems\rev{, impact oscillators, and stick-slip systems}{}.
\end{abstract}

\section{Introduction}
\label{sec:orgc48de16}
Experimental testing plays an important role in the development, calibration, and validation of mathematical models of physical systems.
Systematic testing is of particular importance for nonlinear systems, where the presence of bifurcations can result in a significantly different response for a small parameter change.
A range of bifurcations can be seen in nonlinear systems.
These can pose challenges to simple parameter-sweep tests, with issues including a system going unstable, or jumping to different response regimes to those of interest, and not returning due to hysteresis.
Furthermore, with open-loop testing, bifurcation structures and system dynamics can only be inferred from stable responses.
Unstable responses contain much information for parameter fitting \cite{lee2022modelling,thothadri2005nonlinear,lee2022analysis}, model identification, and inferring separatrices; however they are not normally observable.

Control-based continuation (CBC) is able to systematically reveal both stable and unstable dynamics within an experiment, through a combination of feedback control and path-following methods \cite{sieberControlBasedBifurcation2008}.
It is a model-free experimental method, so that the resulting bifurcation analysis does not depend on modelling assumptions or fitted parameters.
CBC uses feedback control to isolate the dynamics of interest, and to stabilise unstable responses.
\rev{}{A selection of advanced controller strategies have been investigated for CBC, including adaptive and model-predictive control~\cite{li2020adaptivea,li2021adaptive,li2022model,de2022control}.}
Appropriate control targets are solved for iteratively, using standard numerical methods and, when necessary, Fourier-Galerkin discretisation.

CBC has seen significant success on mechanical systems, with examples including resonance curves of a harmonically forced nonlinear oscillator \cite{bartonControlbasedContinuationBifurcation2017}, mapping out fold bifurcations in two parameters \cite{rensonExperimentalTrackingLimitpoint2017}, and finding solution branches of self-excited oscillations \cite{lee2020reducedorder}.
\rev{}{It has been applied to systems exhibiting complex dynamics, such as high levels of friction~\cite{kleymanApplicationControlBasedContinuationCharacterization2020}, isolated vibration modes~\cite{kleymanExperimentalApplicationControlBasedContinuation2021}, and the collective dynamics of pedestrian flows~\cite{panagiotopoulosControlbasedContinuationUnstable2017,panagiotopoulos2022continuation}.}
Solution sets such as equilibria, bifurcations, and periodic orbits can be traced out across a parameter range, uncovering the bifurcation structure of a system, even when faced with noise \cite{de2022control,schilderExperimentalBifurcationAnalysis2015}.
As well as tracing out response curves, features such as backbone curves~\cite{renson2016robust} and dynamical stability \cite{bureauExperimentalBifurcationAnalysis2013,bartonControlbasedContinuationBifurcation2017,rensonApplicationControlbasedContinuation2019} can be inferred.
The collected data can be used for advanced system identification and parameter estimation \cite{beregiRobustnessNonlinearParameter2021,song2022bayesian,de2022control}.

Phase-locked loop (PLL) controllers provide an alternative to CBC when the response curves of interest can be parameterised in terms of shifts between forcing and response phase \cite{peter2017excitation,mojrzisch2012experimental}.
\rev{}{PLL methods have found great success in experimentally computing the backbone curves and nonlinear frequency response functions of mechanical systems~\cite{mojrzisch2016phase,peter2016tracking,denis2018identification}.
CBC and PLL methods both provide consistent, agreeing results~\cite{muller2021comparison,abeloos2022consistency}, and can be coupled together in cases where PLLs fail~\cite{abeloos2023experimental}.}
Nevertheless, CBC is more general than PLL \rev{}{testing, and able to compute a wider range of features}, as it does not require a phase-parameterisation of the response curve.

\begin{figure*}[]
\centering
\includegraphics[width=.7\textwidth]{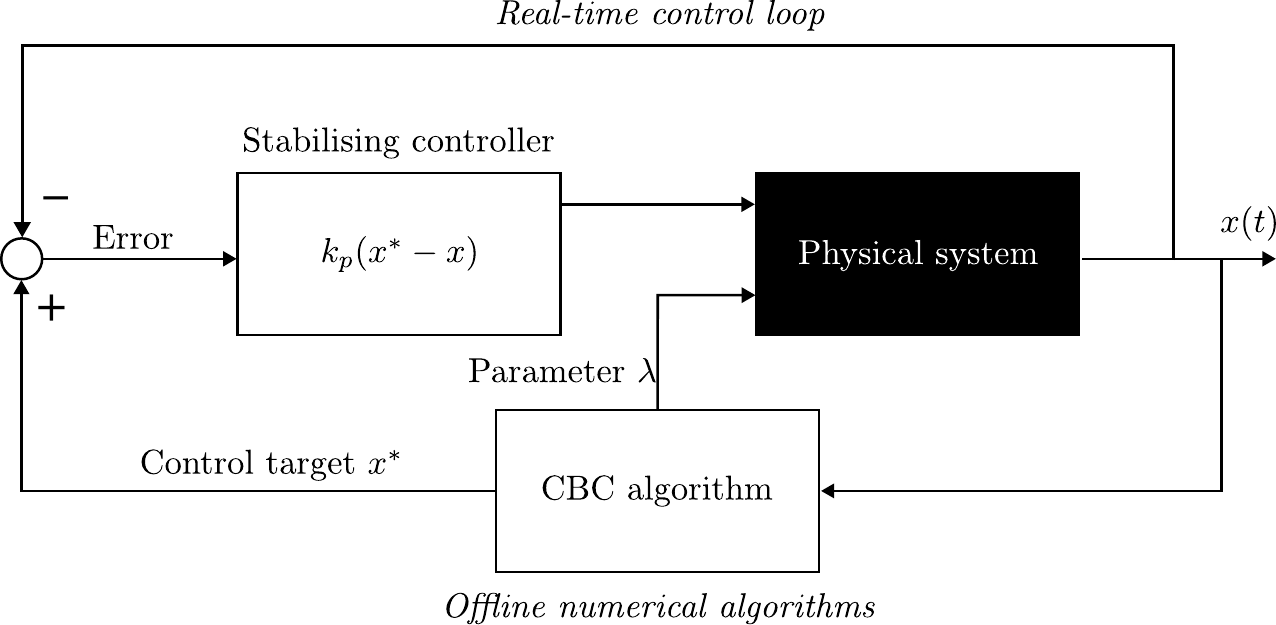}
\caption{\label{fig:org26d902f}Block diagram for control-based continuation experiments. A controllable black-box system of interest is combined with a feedback controller, which provides a stabilising input to the system. The CBC algorithm is used to iteratively update the control target \(x^*\), and the system parameter \(\lambda\).}
\end{figure*}

Slow-fast dynamics appear regularly throughout the physical sciences, whereby state variables evolve at disparate rates.
In engineering, large control gains can cause timescale separations in controlled systems, causing slow-fast responses \cite{kokotovic1999singular}.
In physics, the behaviours of a laser under optical feedback can show mixed-mode oscillations -- an example of multiple-timescale dynamics \cite{marino2011mixedmode}.
Biochemical systems such as gene networks often show a range of expression rates, which produce dynamics over well-separated timescales \cite{chen2002model}.
Similarly, chemical reactions can produce relaxation oscillations, as a result of their slow-fast dynamics; the Oregonator \cite{field1974oscillations} is a classical model of a chemical reaction undergoing relaxation oscillations.
Further examples of systems demonstrating slow-fast dynamics and relaxation oscillations include Josephson junctions in electronics \cite{vernon1968relaxation}; predator-prey systems \cite{liu2003relaxation}; and models of societal growth \cite{milik1996slowfast}.
See \cite{kuehn2015multiple} for additional examples, as well as a detailed coverage of slow-fast dynamical systems.

In principle, CBC is directly applicable to relaxation oscillations, and multiple-timescale systems more generally.
Nevertheless, the tools are limited in practice.
Numerical differentiation is required within the correction step of the continuation.
Due to the model-free nature of CBC, this must be performed using finite differences.
Continuation of periodic responses requires discretisation, with each additional discretisation coefficient imposing an additional system evaluation during finite differencing.
Numerical solutions therefore become impractically slow, when using the large discretisation sizes necessary for relaxation oscillations.
The extra computations also introduce more opportunities for error from noise and system drift \cite{schilderExperimentalBifurcationAnalysis2015}.
Consequently, parsimonious discretisation methods are required to successfully apply CBC to slow-fast systems.

A further issue is that multiple-timescale systems are often autonomous, with no explicit time-dependency.
Oscillatory solutions of autonomous systems comprise of a family of phase-shifted trajectories, making it an \rev{ill-posed}{underdetermined} problem to solve for a single oscillatory response.
Some form of phase condition is needed to extract a unique solution.
Options include the integral phase-condition, used as standard in numerical continuation \cite{kuznetsovElementsAppliedBifurcation2013} and demonstrated numerically with CBC \cite{sieberControlBasedBifurcation2008}, or phase-plane CBC \cite{lee2022analysis}, which is ideally suited for mechanical systems.
However, both methods have limitations -- the integral phase condition can lead to poor numerical stability with CBC, and phase-plane CBC can result in larger discretisation sizes than necessary when applied to relaxation oscillations.

To overcome these problems, we develop an alternative discretisation method and phase constraint for CBC.
B-splines are used in place of the Fourier basis, so that the discretising basis functions can be tailored to the signals of interest.
Optimisation-based techniques are introduced for producing an adaptive discretisation.
The combination of basis functions and adaptation methods allows for an accurate discretisation, with lower dimensionality than Fourier methods.

In addition, we propose an angle-based parameterisation of the discretised control target.
This acts similarly to the recently developed phase-plane CBC \cite{lee2022analysis} to provide a phase constraint for CBC, however our method can lead to lower discretisation sizes on slow-fast signals.
Our angle-based phase constraint does not require knowledge of oscillation frequencies, improving both the speed and accuracy of numerical solutions when compared to the integral phase constraint.
Together, these methods open up possibilities for faster, more accurate CBC experiments on relaxation oscillations.

This work is structured as follows.
Section \ref{sec:orgf215028} introduces the workings of CBC, and the limitations of standard discretisation methods.
Section \ref{sec:org369fcd1} discusses how to construct periodic B-spline models, and section \ref{sec:org0f45276} explains how these may be used as an alternative to the Fourier basis for Galerkin-discretisation.
Section \ref{sec:org9ae8ca1} considers the selection of an optimal B-spline basis for the problem of interest, and how to adapt the discretisation to new data.
Next, section \ref{sec:org10e6e5d} introduces the method of angle-encoded control targets as a phase constraint for autonomous systems. 
Our techniques are demonstrated through simulations on a model of a synthetic gene oscillator and chemical reaction network in section \ref{sec:org6838ea7}.
Section \ref{sec:orgd7cc553} discusses some pertinent aspects of our methods, and section \ref{sec:org0fa60be} concludes the work.

\section{\rev{CBC for}{Control-based continuation of} relaxation oscillations}
\label{sec:orgf215028}

Pseudo-arclength continuation systematically computes the solution manifolds of an underdetermined system \cite{seydel1991tutorial,meijerNumericalBifurcationAnalysis2009}.
With CBC, response families of an experimental system are defined through a specific type of control target, referred to as a noninvasive target, and traced out using pseudo-arclength continuation.
Fig. \ref{fig:org26d902f} depicts a block-diagram of a typical CBC experiment.
It consists of three parts: a controllable system, a stabilising feedback controller operating in real-time, and a suite of numerical methods, operating without time constraints to detect and trace dynamics of interest \cite{sieberControlBasedBifurcation2008,bartonControlbasedContinuationBifurcation2017}.
The controller is used to probe and manipulate the dynamics of the system.
If properly designed, it stabilises any unstable responses, and steers the system towards the dynamics of interest.
A controlled system maps a control target \(x^*(t)\) to an observed output \(x(t)\), referred to as the input-output, or IO-mapping.
Continuation equations are defined on the IO-map.
Parameters and control targets are chosen using a numerical solver within a pseudo-arclength continuation, to locate responses that are intrinsic to the uncontrolled system.

CBC seeks noninvasive control \cite{sieberControlBasedBifurcation2008}.
A controller is noninvasive if it stabilises a limit cycle or equilibrium that exists within the uncontrolled system, without changing its position in parameter space.
Hence, for control target \(x^*(t)\) and measured system response \(x(t)\), a sufficient condition for an applied control action \(u(x^*, x, t)\) to be noninvasive is if \(x^*\) is stabilised, and \(u\equiv 0\) for all time.
Characteristics such as the geometry, oscillatory period, and location of steady-state system responses remain unchanged in the presence of noninvasive control.
However, unstable features are stabilised, and become directly observable.

Fixed points of the IO-map are any control targets that are exactly tracked by the system output, and therefore satisfy \(x^* = x\).
For a proportional or proportional-plus-derivative controller, the total control action \(\| u(x^*, x, t)\|\) is zero if and only if \(x^*(t) = x(t)\), so that noninvasive control is given by fixed points of the IO-map.
Hence, noninvasive control targets are found by using a nonlinear solver to solve for \(x(t) - x^*(t) = 0\), guaranteeing zero control input and natural system dynamics.

\begin{figure*}[th!]
\centering
\includegraphics[width=.9\linewidth]{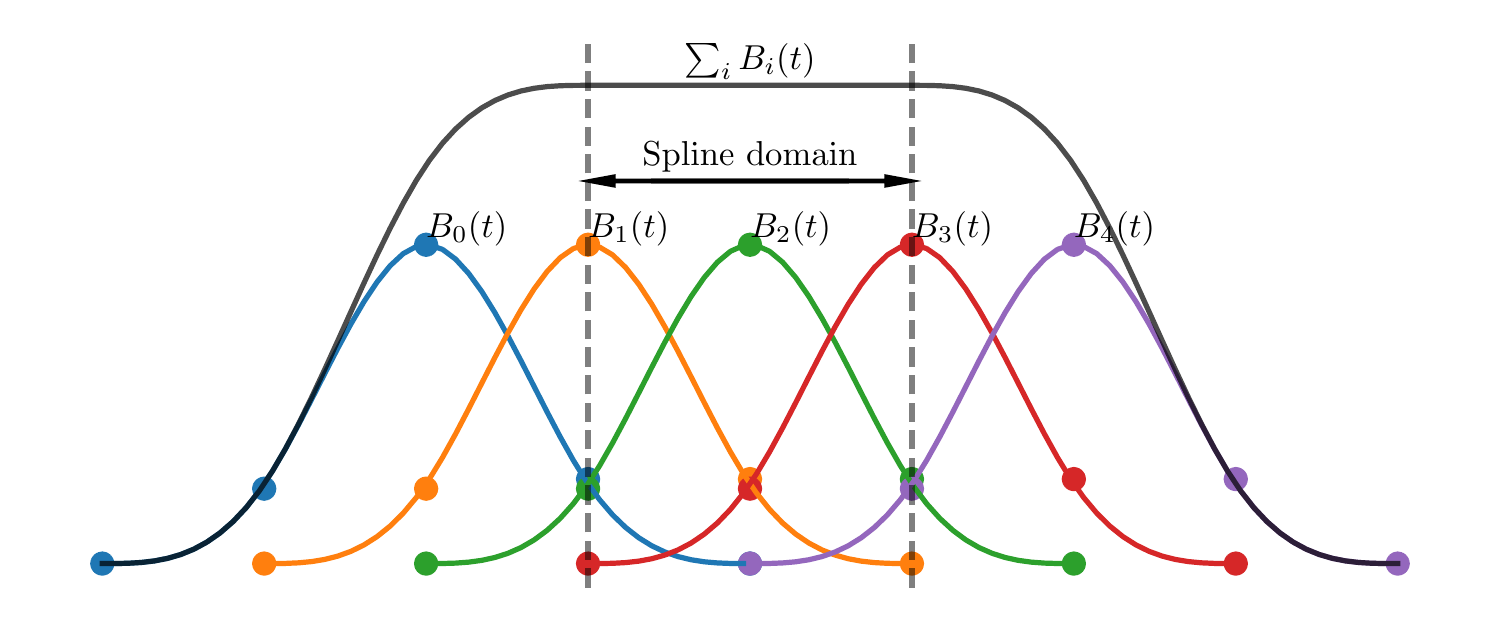}
\caption{\label{fig:orgae46232}\rev{A cubic spline model with equispaced knots. Observe that the B-spline functions extend beyond the spline domain. The cubic spline curve is only valid when there are exactly four active B-spline functions.}{An illustrative cubic spline model, constructed from five basis functions, each with equispaced knots. Each cubic B-spline is constructed from five consecutive knots. The range of these five knots gives the domain of that B-spline function. Within this domain, the spline is referred to as being active. The domain of the cubic spline curve itself, shown here in dotted grey, is the region within which exactly four B-spline functions are active. As a result, knots must be placed on or outside the domain of the spline model, to provide the appropriate number of B-splines within the domain. These are referred to as exterior knots.}}
\end{figure*}

Solving for fixed-points is performed using standard numerical methods, such as Newton or quasi-Newton iterations \cite{sieberControlBasedBifurcation2008}.
These cannot be applied directly when \(x^*\) is oscillatory; instead, the zero-problem must be discretised.
All current CBC implementations use a Fourier discretisation, in which control targets and system responses are represented by truncated Fourier series, and equality is sought between Fourier coefficients \cite{sieberControlBasedBifurcation2008}.
Assume the control target and measured system output are well-approximated by their first \(n\) Fourier modes.
Expanding both functions in the Fourier basis gives
\begin{align}
x^*(t) &= a_0^* + \sum_{k=1}^n a_k^*\cos(k\omega t) + b_k^*\sin(k\omega t)~,\\
x(t) &= a_0 + \sum_{k=1}^n a_k\cos(k\omega t) + b_k\sin(k\omega t)~,
\end{align}
where \(\omega\) is the oscillatory frequency.
Then, \(x^*-x=0\) is satisfied when \(a_0^* - a_0 = 0\), \(a_k^* - a_k = 0\) and \(b_k^* - b_k=0\), for \(k\in\{1, \dots, n\}\).
Unlike the original problem, the discretised problem can be solved using standard numerical methods.
Furthermore, the projection onto basis functions provides a degree of noise-averaging, improving robustness with noise-corrupted measurements.

Fourier discretisation has been used successfully in existing CBC experiments, where signals are well-approximated by few Fourier harmonics.
Nevertheless, it presents challenges for systems which undergo rapid changes, such as relaxation oscillations or impacting and friction dynamics.
Such systems produce highly changeable signals, and many Fourier modes are needed for producing an accurate discretisation.
This results in a high-dimensional continuation problem, which in turn slows down experiments and introduces more opportunities for error.
Here, we develop a B-spline discretisation scheme to overcome this issue.

\subsection{Periodic B-splines}
\label{sec:org369fcd1}
\begin{figure*}[th!]
\centering
\includegraphics[width=.8\textwidth]{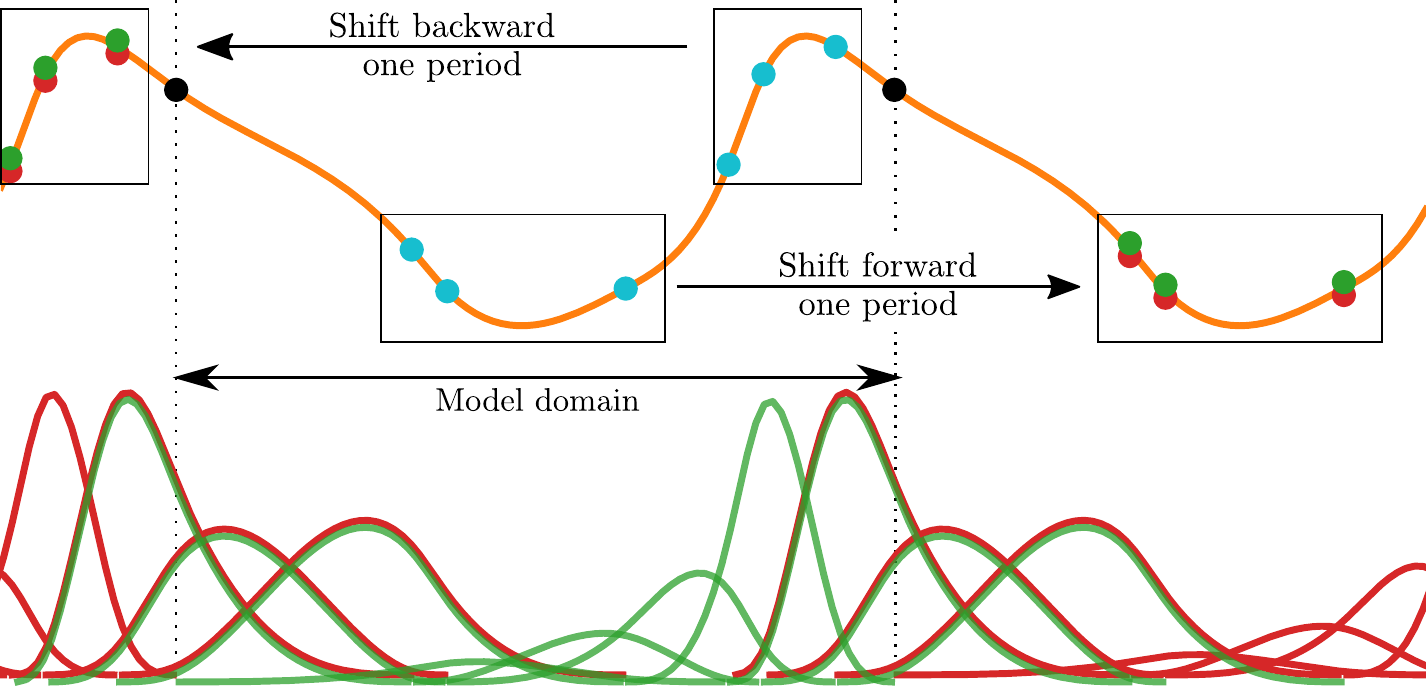}
\caption{\label{fig:org077912f}Knot positions (dots) and basis functions (green, red curves) of a B-spline model with periodic boundary conditions (orange curve). The last three interior knots of the spline model (blue dots) are shifted back by one period to get the left exterior knots (green dots), and first three interior knots are shifted forward by one period to get the right exterior knots. Boundary knots (black dots) mark the points at which the resulting periodic boundary conditions are applied. Note that if an identical spline model were to be created one period in the past or future, the interior knots of the adjacent periods (red dots) would overlap exactly with the interior knots of this period; pairs of boundary knots would overlap exactly. As a result, the basis functions of the adjacent periods (red curves) would exactly overlap those of this period (green curves), providing the equalities necessary for periodic boundary conditions.}
\end{figure*}

A spline is a maximally smooth, piecewise-polynomial curve.
Splines are defined by choosing a set of knot points and connecting them together with sections of polynomial.
Polynomial coefficients are determined partially by ensuring that each polynomial starts and ends at its boundary knots, and meets its neighbouring polynomials smoothly.
In addition, a pair of boundary conditions are required.
In this work we assume periodic boundary conditions.

Spline curves are often expressed as a weighted sum of basis functions, to simplify computations.
B-splines are a family of minimal-support basis functions for spline curves \cite{deboorPracticalGuideSplines1980}.
The B-spline basis is constructed from a set of scalar-valued knots, defining the \(x\) positions at which each section of polynomial meets.
The \(y\) positions are determined by the basis function coefficients.

The B-spline basis can be obtained through the following recurrence relation (see \cite{deboorPracticalGuideSplines1980} for a detailed discussion and proof).
Let \(\Xi = \{\xi_0 \leq \dots \leq \xi_n\}\) be an ordered set of knots.
Given \(\Xi\), the first-order B-splines, consisting of piecewise-constant terms, are defined as
\[B_{i,1}(t) =
\begin{cases}
1, &\quad \text{if } \xi_i \leq t \leq \xi_{i+1} \\
0, &\quad \text{otherwise.}
\end{cases}\]
Higher-order spline functions are constructed recursively; let
\[
\omega_{i,k}(t) = \begin{cases}
  \frac{t - \xi_i}{\xi_{i+k-1} - \xi_i}, &\quad t_{i+k-1} \neq t_i \\
  0, &\quad \text{otherwise.}
\end{cases}\]
Then the \(i\)th B-spline of order \(k\), defined on knots \(\Xi\), is given as
\[
B_{i,k}(t) = \omega_{i,k}(t)B_{i,k-1}(t) + \big[1-\omega_{i+1,k}(t)\big]B_{i+1,k-1}(t)~.
\]
B-spline functions are generally evaluated using the Cox-de Boor algorithm \cite{coxNumericalEvaluationBsplines1972}.
Implementations are provided as standard in many scientific computing packages.

Fig. \ref{fig:orgae46232} depicts a cubic B-spline model with equidistant knots.
\rev{}{By definition, each B-spline function has minimal support across a domain determined by its knots.
Within the domain, the spline is referred to as being active.}
It is observed that exactly four B-splines are active at any point within the domain of the spline curve for which the B-splines form a basis.
A consequence of this is that knots are also required outside of the spline domain, to ensure the correct number of active B-splines within the domain.
We refer to the knots outside of the domain as exterior knots, and knots within the domain as interior knots.
Exterior knots dictate the boundary conditions of the spline model.
Here we recall how exterior knots can be chosen to give periodic boundary conditions, using the method of periodic extension \cite{deboorPracticalGuideSplines1980}.
Section \ref{sec:org9ae8ca1} then considers how to select optimal interior knots.

As depicted in Fig. \ref{fig:org077912f}, periodic boundary conditions require the exterior knots and B-spline functions from one period to overlap with those of adjacent periods.
This is achieved simply by shifting appropriate sets of interior knots by one period \cite{deboorPracticalGuideSplines1980}.
Given a knot on the lower and upper domain boundary, and three or more interior knots, periodic exterior knots for a cubic set are found by
\begin{itemize}
\item shifting the last three interior knots (blue dots in Fig. \ref{fig:org077912f}) backwards by one period, to produce the first three exterior knots (green dots);
\item shifting the first three interior knots (blue dots) forward by one period, to produce the last three exterior knots (green dots);
\item enforcing equality between respective pairs of coefficients for the first and last three B-splines.
\end{itemize}
This ensures equality between the B-spline functions at the start and end of the period, and produces periodic boundary conditions.

More rigorously, consider a cubic spline curve
\begin{equation}
b(t) = \sum_{i=-3}^n \beta_i B_i(t),
\label{eq:thisb}
\end{equation}
defined on the domain \(t\in[0,T]\), with periodic boundary conditions.
The boundary conditions require \(b(0) = b(T)\), and likewise for all available derivatives.
Exactly four basis functions are active at the start of the domain, and equivalently at the end.
Considering only these active basis functions, it is seen that periodic boundary conditions require
\begin{equation}
\sum_{i=-3}^{0}\beta_iB_i(t+T) = \sum_{i=n-3}^n \beta_iB_i(t)~,
\end{equation}
implying
\begin{equation}
\beta_{-3} = \beta_{n-3} , \quad \dots, \quad \beta_0 = \beta_n
\end{equation}
and
\begin{equation}
B_{-3}(t+T) = B_{n-3}(t), \quad \dots, \quad B_0(t+T) = B_n(t).
\end{equation}
Each cubic B-spline is defined over a set of five knots; let B-spline \(B_i(t)\) be constructed from knots \(\Xi^i = \{\xi^i_1 \leq \xi^i_2 \leq \xi^i_3 \leq \xi^i_4 \leq \xi^i_5\}\).
Shifting a set of B-spline knots by some amount \(T\) produces an equivalent shift in the resulting B-spline function.
That is, given \(\Xi_T^i = \{\xi^i_1 + T\leq \xi^i_2 + T \leq \xi^i_3 +T \leq \xi^i_4 +T \leq \xi^i_5 + T\}\), the corresponding B-spline function \(B_i^T(t)\) satisfies \(B_i^T(t+T) = B_i(t)\).
Hence, for \(B_{-3}(t+T) = B_{n-3}(t)\) to be satisfied, we must have \[ \Xi_{n-3} = \Xi_{-3} + T , \quad \Xi_{n-2} = \Xi_{-2} + T, \quad \text{etc.} \]

\subsection{B-spline discretisation for CBC}
\label{sec:org0f45276}
B-spline discretisation for CBC proceeds in much the same way as Fourier discretisation.
Noninvasiveness is sought by solving for equality between the control target and system output, when proportional or proportional-plus-derivative control are used.
This is achieved by projecting a system response onto a B-spline basis, and seeking equality between the target and response coefficients.
Control targets are represented as a sum of basis functions \(B_i\) weighted by coefficients \(\beta_i\), as given by
\begin{equation}
x^*(t) = \sum_{i=1}^n \beta_i^* B_i(t)~.
\label{eq:simplespline}
\end{equation}
Periodicity of the control target is guaranteed by selecting basis functions \(B_i\) for periodic boundary conditions, as discussed in section \ref{sec:org369fcd1}.
Target \(x^*(t)\) is used to control the system of interest, using a proportional or proportional-plus-derivative strategy.
Recorded samples \((t_j, x_j)\) from the system output \(x(t)\) are projected back onto the basis functions \(B_i\) using least-squares.
Given samples \((t_j, x_j)\) for \(j\in\{1,\dots, N\}\), let
\begin{align}
[\mathcal{B}]_{u,v} &= B_v(t_u) ~, \\
\beta &= [\beta_1, ~ \beta_2, ~ \dots, ~ \beta_n]^T ~, \\
X &= [x_1, x_2, \dots, x_N]^T ~.
\end{align}
Equation \eqref{eq:simplespline} suggests \(X = \mathcal{B}\beta\), however given sufficient data, this is overdetermined, with a least-squares error \(e\) of
\begin{equation}
e = \| X - \mathcal{B}\beta \| ^2~.
\label{eq:lsqerror}
\end{equation}
Coefficients \(\beta_\mathrm{lsq}\) that minimise this error are given by the least-squares Eq.
\begin{equation}
\beta_\mathrm{lsq} = \left(\mathcal{B}^T \mathcal{B}\right)^{-1} \mathcal{B}^T X~.
\label{eq:lsqsoln}
\end{equation}
Hence, the B-spline representation \(x(t)\) of samples \((t_j, x_j)\) is given by
\begin{equation}
x(t) = \sum_{i=1}^n \beta_i B_i(t)~,
\end{equation}
for \(\beta_\mathrm{lsq} = [\beta_1,~ \dots,~ \beta_n]\).
Noninvasiveness is met when \(\beta_i^* - \beta_i =0\) for all \(i\in\{1,\dots,n\}\).

Here we exclusively consider cubic spline models.
Nevertheless, the methods presented in this work can be implemented using splines of any order.

\subsection{Selecting and adapting interior knots}
\label{sec:org9ae8ca1}
\begin{algorithm*}[t]
    \begin{algorithmic}
	\Function{Objective}{interior knots, reference signal}
	    \State $\Xi \gets $ full knot set \Comment{Get boundary and exterior knots}
	    \State Get B-spline basis from $\Xi$ \Comment{Eg. Cox-de Boor algorithm}
	    \State $e_\Xi \gets \|X - \left(\mathcal{B}^T \mathcal{B}\right)^{-1} \mathcal{B}^T X \beta \|^2$ \Comment{Least-squares error, from Eq.s \eqref{eq:lsqerror} and \eqref{eq:lsqsoln}}.
	    \State \Return $e_\Xi$
	\EndFunction
    \end{algorithmic}
    \begin{algorithmic}
        \Procedure{Knot placement}{$n$, reference signal, repeats}
	    \State $e_\text{min} \gets \infty$
	    \State $s \gets 1$
	    \Repeat
	        \State $\Xi_0 \gets \mathcal{U}_{[0,1]}^n$
		\Comment{Draw initial knots at random from \(n\)-dimensional uniform distribution}
		\State $\Xi_\text{opt} \gets \argmin_\Xi$ \Call{Objective}{$\Xi, \text{reference signal}$}
		\Comment{Numerically optimise $\Xi$ using L-BFGS-B, with $\Xi_0$ as initial guess}
		\State $e_\Xi \gets $ \Call{Objective}{$\Xi_\text{opt}, \text{reference signal}$}
		\If{$e_\Xi < e_\text{min}$} \Comment{Store best-performing knots}
		    \State $e_\text{min} \gets e_\Xi$
		    \State $\Xi_\text{best} \gets \Xi_\text{opt}$
		\EndIf
	        \State $s \gets s+1$
	    \Until{$s > \text{repeats}$}
	\EndProcedure
    \end{algorithmic}
\caption{Pseudo-code for selecting spline knots to minimise least-squares error.}
\label{alg:optiknots}
\end{algorithm*}

B-spline knots act analogously to a mesh in finite element methods.
While uniformly spaced knots could be used, it is typically beneficial to tailor the knots to the problem of interest.
Smaller discretisation sizes are almost always achieved when knot positions are selected specifically for the system of interest.
Here we discuss a method to select an optimal set of spline knots, and to adaptively update them throughout a continuation.

Many methods have been proposed for selecting spline knots.
Free-knot methods choose knots either through one of many optimisation-based techniques \cite{valenzuelaEvolutionaryComputationOptimal2013,kangKnotCalculationSpline2015,schwetlickLeastSquaresApproximation1995,juppApproximationDataSplines1978}, using various heuristics \cite{liHeuristicKnotPlacement2004,michelNewDeterministicHeuristic2021}, or by creating probabilistic distributions over knot sets \cite{dimatteoBayesianCurvefittingFreeknot2001,lindstromBayesianEstimationFreeknot2002,mamicAutomaticBayesianKnot2001}.
For CBC, the experimenter should be free to determine their desired discretisation size.
As such, optimisation-based knot placement strategies are appealing, whereby the experimenter chooses the number of knots, and the experimental results inform where they are best placed.
We use numerical optimisation to automatically place a predetermined number of knots, in order to minimise the least-squares error between a set of reference data, and the best-fit spline model of those data.

The optimisation procedure is summarised in algorithm \ref{alg:optiknots}, and knot selection proceeds as follows.
For each iteration of a local optimiser, a set of interior knots are produced.
From these, exterior and boundary knots are added, and used to construct a set of basis functions.
A least-squares fit to a reference signal is produced using these basis functions, and the local-optimiser algorithm performs further iterations to find the interior knots that minimise this fitting error.
We choose a bounded limited-memory BFGS (L-BFGS-B) local optimiser, implemented in SciPy \cite{virtanenSciPyFundamentalAlgorithms2020}.
To avoid local minima, the optimisation is repeated several times from random initial conditions, and the best result is selected as the final knot set.
The reference signal can be measurements from an uncontrolled system evaluation when initialising knots, or the last accepted continuation solution when following the knot adaptation procedure discussed next.
As observed in Fig. \ref{fig:org077912f}, optimal knots tend to be placed around turning-points within a signal.

To ensure the optimiser proposes valid knots, constraints are applied to the problem, to restrict proposed knots to within the model domain.
Samples are projected onto a single period.
We choose, without loss of generality, to rescale that period onto the unit interval.
Each B-spline knot \(\xi\) can then be constrained to \(\xi\in[0,1]\), providing bounds on the search space.

The objective function requires no additional evaluations of the controlled or uncontrolled system.
Optimisation can be performed entirely offline in a matter of seconds.
Globally optimal knots are typically found.
An alternative to random restarts is a global optimiser such as simulated annealing.
In our tested cases, this is found to generally produce the same knots as random restarts, however random restarts produced satisfactory results more rapidly for larger discretisation sizes.

Continuation solutions will change during an experiment, as parameters are varied.
It is preferable to adapt the discretisation basis throughout the continuation.
We achieve this by running an optimisation step on the current knots every time a new solution has been accepted.
We choose either a sequential least-squares quadratic programming (SLSQP), or constrained trust-region optimisation method for this, both implemented in SciPy \cite{virtanenSciPyFundamentalAlgorithms2020}.
A new, different set of globally optimal knots may emerge during an experiment.
The experimenter may therefore choose to perform a full multi-restart knot optimisation when the discretisation error increases beyond a predetermined threshold.
When using adaptive discretisation, care must be taken to ensure \rev{a consistent set of basis functions are used within any given computation.}{that the same basis functions are used to discretise all solutions within a given prediction-correction step; pairs of previous results cannot be used for secant prediction if their discretised solutions are obtained from different basis functions.}
Consistency in secant predictions is obtained by rediscretising previous solution data, and projecting all samples onto the newly updated basis functions.
Hence, an accepted solution is used as a reference signal for selecting knots, then the same set of basis functions are used throughout all calculations in the next prediction-correction step.

\subsection{Angle-encoding phase constraint}
\label{sec:org10e6e5d}
\begin{figure*}[]
\centering
\includegraphics[width=\textwidth]{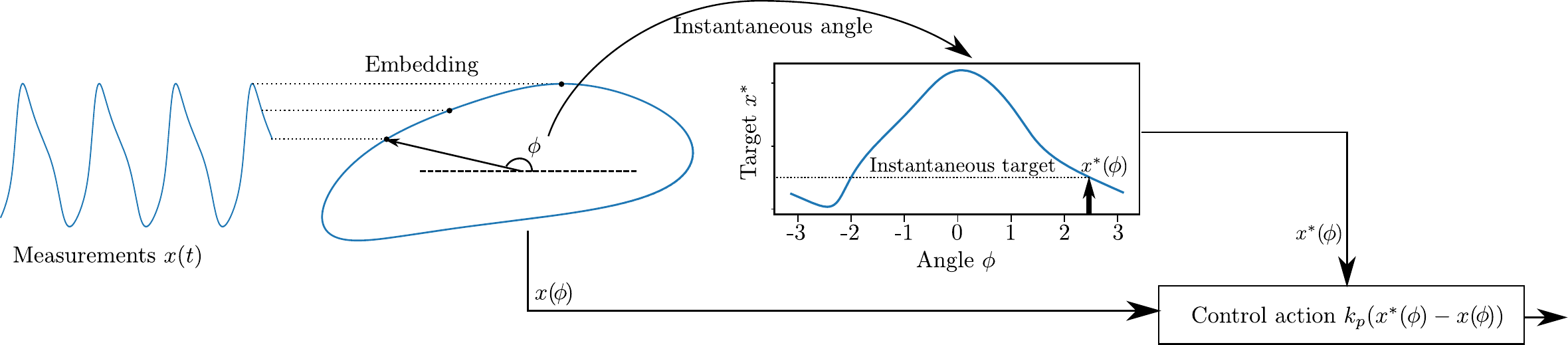}
\caption{\label{fig:org0691e14}Schematic diagram of the angle-encoded control strategy. Measurements \(x(t)\) from the controlled system are embedded in a planar limit cycle. The angle of the current embedded state is used to find the instantaneous control target. The control target is sent to a feedback controller, compared to measurements, and fed back into the system as a control action.}
\end{figure*}

A phase-constraint is required to produce a unique solution for the continuation of oscillations in an autonomous system.
Here, we obtain a unique solution by indexing control targets in terms of an angle-based independent variable, instead of time.
As CBC solutions now lack any time-dependency, phase shifts cease to be meaningful, and oscillatory solutions become (locally) unique.
We refer to our method as angle-encoding, since control targets are solved for in the angle-domain.
Noninvasiveness is achieved when equality is reached between the angle-encoded system response and control target.

Angle-encoding proceeds by taking time-dependent samples \((t_i, x_i)\), and converting them to time-independent samples \((\phi_i, x_i)\).
This is performed by replacing the time index variable \(t\) with an angle-based independent variable \(\phi\).
Control targets are then constructed to give the desired system output at some state-angle \(\phi\).

The full angle-encoding procedure is as follows.
Some system response \(x(t)\) is measured.
A second variable \(z(t)\) is also required\rev{; this is}{,} referred to as the embedding variable.
The embedding variable is used to reconstruct a planar limit cycle from the observed data; it can come from an explicitly measured state variable\rev{}{ in an experiment}, or from \rev{processing of measurements \(x(t)\).}{a proxy for the state, such as delay or derivative coordinates.}
Embedding variable \(z(t)\) is used in conjunction with state variable \(x\) to calculate an instantaneous state angle \(\phi\).
\rev{}{The chosen embedding must produce a unique mapping from angles to control targets.}
\rev{Here, we use a second explicitly observed state variable \(y\) to reconstruct a planar limit cycle.}{}

For some embedding variable \(z(x,y)\), the angle of an embedded state \((x(t), z(t))\) is given by
\begin{equation}
\phi(x, z) = \mathrm{atan2}\left(\sigma\left(z - \mu_z\right), x - \mu_x\right).
\label{eq:atan}
\end{equation}
Parameters \(\mu_x\) and \(\mu_z\) shift the embedded limit cycle to encircle the origin, and \(\sigma\) scales it.
Origin choice \((\mu_x, \mu_z)\) is found to affect the performance of the controller, as discussed further in section \ref{sec:orgd7cc553}.
Scale-factor \(\sigma\) is beneficial in the cases where \(x\) and \(z\) have significantly different amplitudes.
Such cases often arise with slow-fast systems, when using derivative coordinates as an embedding scheme.
The scale factor ensures angles are well-distributed across the interval \([0, 2\pi)\), so that the angle-encoded control target can be represented in an effective manner.
A good choice for \(\sigma\) is anything that brings \(x(t)\) and \(\sigma z(x, y)\) into comparable amplitudes.
Principal component analysis may be used in place of scale factor \(\sigma\), to transform embedded states when \(x\) and \(z\) are highly correlated, for example in delay embeddings.

Angle-encoded signals are time-independent, so phase-shifts are no longer meaningful.
Solution phases are instead determined by the combination of system dynamics and angle definition.
As such, angle-encoding ensures that solving for an oscillatory response is a well-posed problem, with a unique solution.
Furthermore, no knowledge of oscillatory period is required.
This simplifies the continuation equations, by avoiding the need to calculate the \rev{natural}{response} frequency of the system.

\subsection{Angle-based control}
\label{sec:org05ce2db}

\begin{figure*}[th!]
\centering
\includegraphics[width=.8\textwidth]{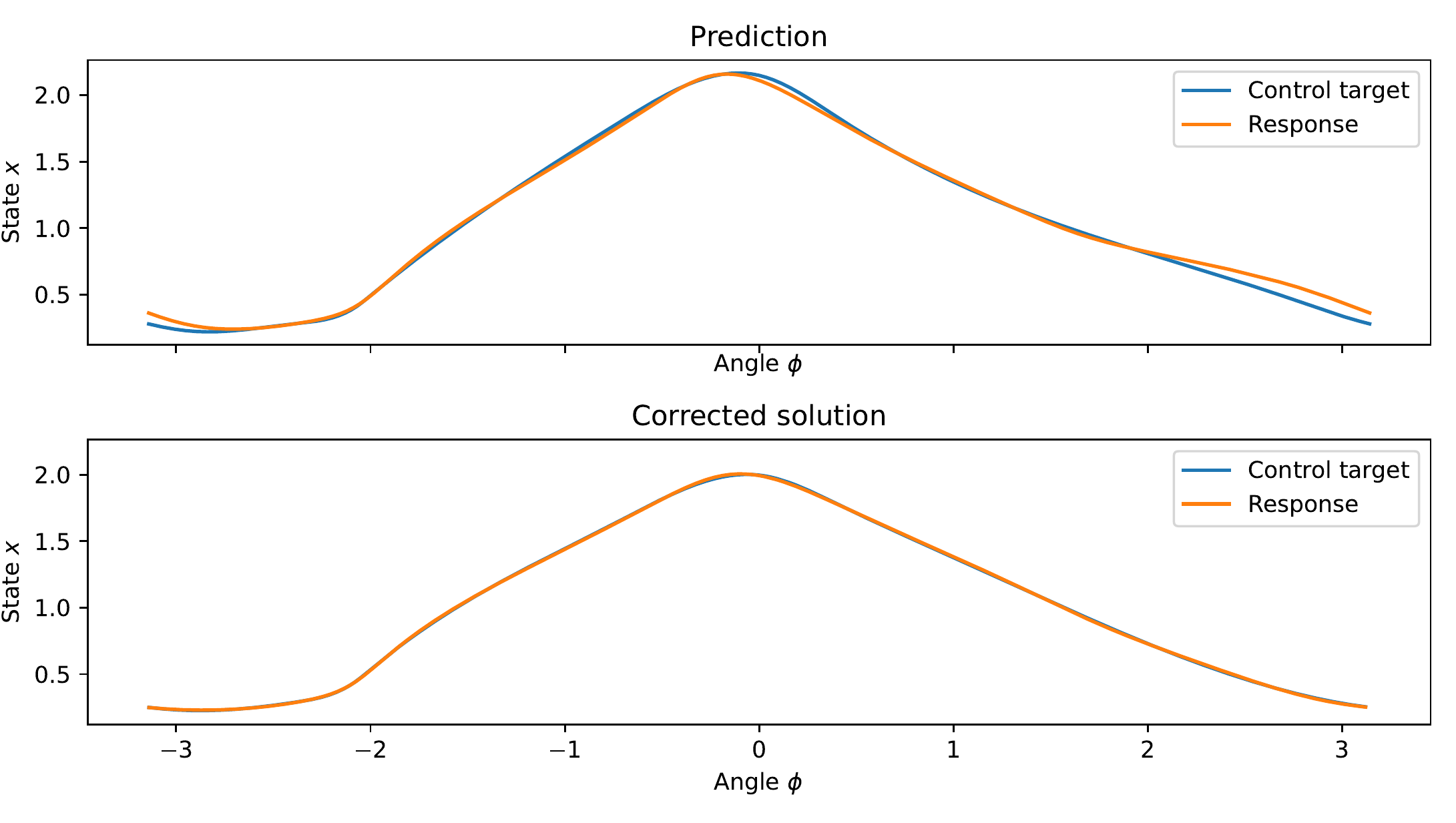}
\caption{\label{fig:org55d0e07}Examples of an angle-encoded control target and system response, before and after a correction step. Noninvasive control is achieved when equality is found between the angle-encoded discretisations of the control target and system response.}
\end{figure*}

Control proceeds as sketched in Fig. \ref{fig:org0691e14}.
At any given instant, the system is observed.
Observed measurements \(x(t)\) are combined with some embedding variable \(z\) and transformed, to reconstruct a planar limit cycle.
The angle of the current state on this limit cycle is calculated with Eq. \eqref{eq:atan}, and a proportional control law of form
\begin{equation}
u(x, z) = k_p\left(x^*(\phi) - x(\phi)\right)
\end{equation}
is chosen.
Control force \(u\) is then fed back into the system, to stabilise target \(x^*\).

The continuation zero problem is replaced by its angle-encoded counterpart, as
\begin{equation}
x(\phi) = x^*(\phi)~.
\label{eq:noninvasiveangle}
\end{equation}
This requires the measured system output to be translated into an
angle-encoded output, in the same way as is performed by the
controller. Eq. \eqref{eq:noninvasiveangle} is discretised and
solved identically to its time-dependent counterpart.

Angle-encoded noninvasiveness behaves exactly the same as temporal noninvasiveness, by guaranteeing zero control action.
An example is plotted in Fig. \ref{fig:org55d0e07}, showing angle-encoded targets and responses.
These are taken from the first prediction-correction step of the gene oscillator CBC simulation discussed in section \ref{sec:orgef114e2}, except with a larger stepsize of 0.75.
Angle-encoded control targets and system responses are shown, both for an uncorrected prediction, and for the same solution after Newton-correction.
It is seen that equality has been found between the angle-encoded control target and system response, giving a noninvasive solution.

\section{Case studies}
\label{sec:org6838ea7}
Our methods are demonstrated through simulations of two slow-fast models -- a gene oscillator, and an oscillating chemical reaction.
Notes on possible experimental implementations are provided, the CBC method is discussed, and simulated results are presented for each system.

The first test-system is taken from synthetic biology.
Mathematical modelling and systems identification are powerful tools within systems and synthetic biology \cite{marucciNanogDynamicsMouse2017,dibernardoPredictingSyntheticGene2012,marucci2011derivation}.
Recent developments in the field have produced a range of novel methods for the automatic feedback control of gene expression and signalling pathways in live cells \cite{ren2017population,postiglione2018regulation,shannon2020vivo,pedone2019tunable,menolascina2014vivo,khazim2021microfluidic,pedone2021cheetah,lugagne2017balancing}.
Accordingly, methods are being actively developed to exploit these control methods, for applying CBC to synthetic gene networks \cite{de2022control}, which could be of great benefit when it is difficult to derive or calibrate models for the biological system of interest.
Here we show numerical results of a CBC experiment, to track oscillations in a modelled slow-fast gene regulatory network.
Results are obtained from a nondimensional model taken from \cite{hasty2002synthetic}.
For proteins \emph{X} and \emph{Y}, with respective (nondimensionalised) concentrations \(x\) and \(y\), the network dynamics \rev{are}{is} given by
\begin{align}
\label{eq:hasty}
\frac{\mathrm{d} x}{\mathrm{d} t} &= \frac{1 + x^2 + \alpha \sigma x^4}{(1 + x^2 + \sigma x^4)(1 + y^4)} - \gamma_x x ~ ,\\
\label{eq:hasty2}
\tau_y \frac{\mathrm{d} y}{\mathrm{d} t} &= \frac{1 + x^2 + \alpha \sigma x^4}{(1 + x^2 + \sigma x^4)(1 + y^4)} - \gamma_y y ~,
\end{align}
with \(\tau_y=10\), \(\alpha=11\), \(\gamma_x=0.105\), \(\sigma=2\); we consider continuation parameter \(\gamma_y \in [0.01, 0.05]\).
All parameter values are as given in \cite{hasty2002synthetic}, except for design parameter \(\tau_y\).
This quantifies the system timescale separation in the system, which we increase to 10 to further widen the timescale separation.
Bifurcation parameter \(\gamma_y\) quantifies the degradation rate of protein \emph{Y}, which can be modified by varying the concentration of isotropyl-\(\beta\)-D-thigalactopyranoside (IPTG) \cite{hasty2002synthetic}.
Recent methods have been developed to modify IPTG concentrations online for cell cultures in microfluidic chambers.
In these cases, IPTG is used as a control input to cells \cite{shannon2020vivo,khazim2021microfluidic}; the same method can be used to change IPTG concentrations for modifying bifurcation parameter \(\gamma_y\).
Existing genomic control experiments use bang-bang control strategies\rev{}{~\cite{lasalle1959time}}; it remains to be seen whether these methods have the necessary fidelity to stabilise unstable oscillations.
Optogenetic methods provide an alternative approach for implementing control strategies, by using various wavelengths of light to control gene expression \cite{tabor2011multichromatic,multamakioptogenetic,chia2022optogenetic}.

Our second test system is the reduced Oregonator model, describing the dynamics of the oscillating, autocatalytic Belousov-Zhabotinsky reaction.
Forcing and feedback control have been used for a variety of experiments on the dynamics of Belousov-Zhabotinsky reactions.
A wide range of feedback controllers have been proposed \cite{fatoorehchi2015feedback}, with common experimental strategies include combining video measurements and optical feedback \cite{kheowan2001spiral,grill1995feedback,goldschmidt1998transition,grill1996spiral,tung2002dynamics}, and by modifying reactant inflow rates \cite{fatoorehchi2015chaos,petrov1993controlling}.
Other experiments use electrical stimulation to provide perturbations into the experiment \cite{showalter1979detailed,kuze2019chemical,peralta2006controlled}.
We consider feedback control of an ODE model, representing a well-stirred Belousov-Zhabotinsky reaction which could be controlled through modification of reactant inflow rates.
Dynamics is governed by a two-variable nondimensionalised Oregonator, obtained through a quasi steady-state approximation of the full system \cite{tyson1980target}.
The model is given by
\begin{align}
\label{eq:oregon}
\varepsilon \frac{\mathrm{d} x }{\mathrm{d} t} &= x(1-x) - f\frac{y(x-q)}{x+q} ~,\\
\frac{\mathrm{d} y }{\mathrm{d} t} &= x - y~.
\end{align}
We take \(\varepsilon=0.1\), \(q=0.025\), and \(f\) is to be varied, as the continuation parameter.

\subsection{CBC parameters}
\label{sec:orgf6bbda2}
\begin{figure*}[th!]
\centering
\includegraphics[width=.8\textwidth]{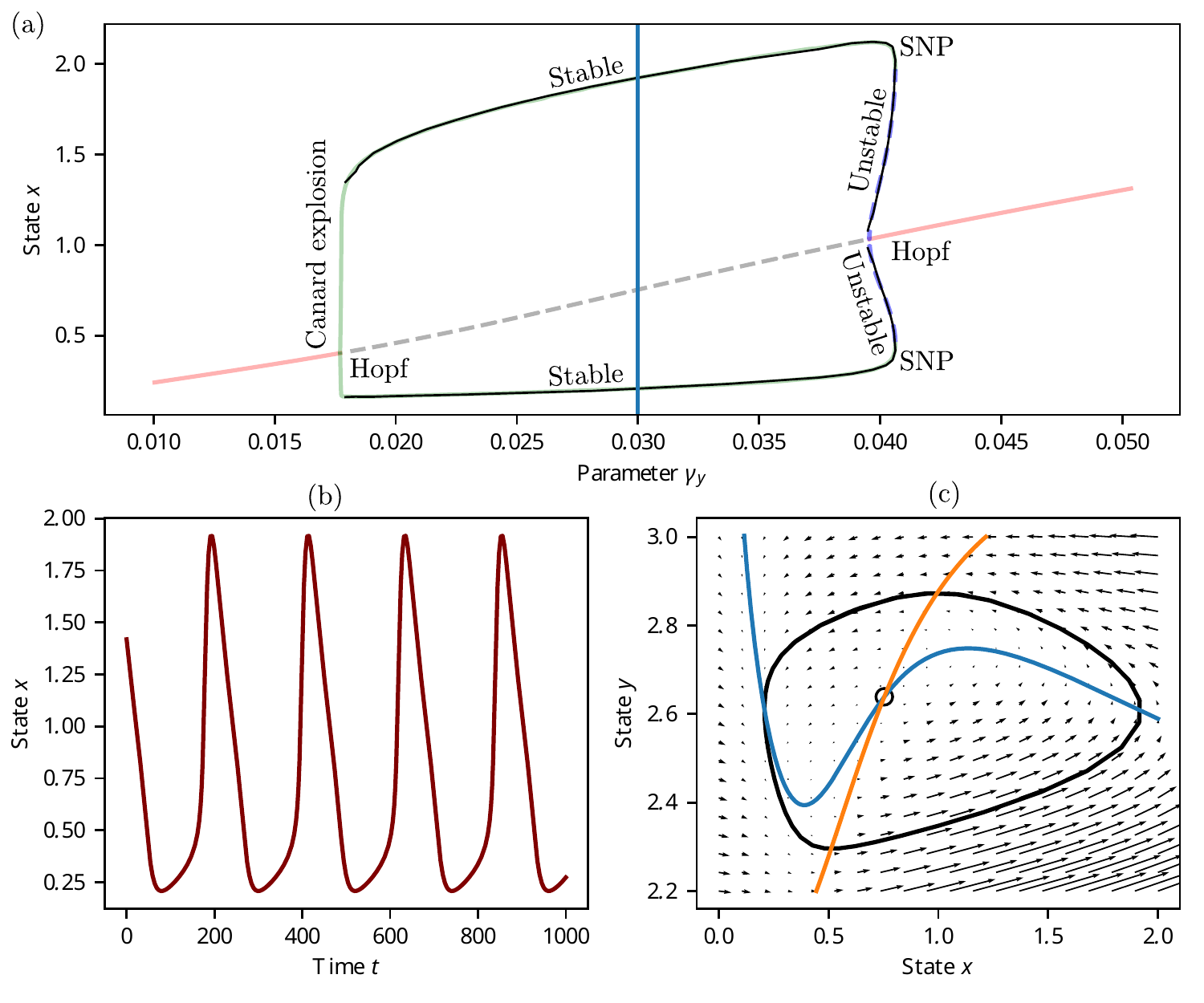}
\caption{\label{fig:org8d4d85c}Simulated CBC results for an oscillating synthetic gene network \cite{hasty2002synthetic}. CBC results are shown on panel (a) in black, and are overlaid on a bifurcation diagram generated in XPPAUTO (version 8.0), showing stable (\rev{red}{pink}) and unstable (\rev{black}{grey}-dotted) equilibria, and stable (green) and unstable (\rev{blue dotted}{purple dotted}) limit cycles. The dynamics of the oscillator are shown in panels (b) and (c), as a time-series plot of the \(x\) variable (panel b), and a phase-plane diagram illustrating the nullclines (blue, orange lines), limit cycle (black line), equilibrium (open circle), and flow-field (black arrows), panel (c). Time series and phase plane plots are produced at \(\gamma_y = 0.03\), indicated by the vertical blue line in the bifurcation diagram. Stable and unstable branches of periodic orbits are marked as such; CBC stabilises the unstable branches, to render them observable. In the open-loop system, branches of periodics emerge at the marked Hopf bifurcations; branch stability switches at the saddle-node of periodic orbits bifurcation (SNP), and on the canard explosion. CBC results are seen to match the XPPAUTO results almost perfectly, with the exception of the canard explosion which CBC cannot follow.}
\end{figure*}

\begin{figure*}[th!]
\centering
\includegraphics[width=.8\textwidth]{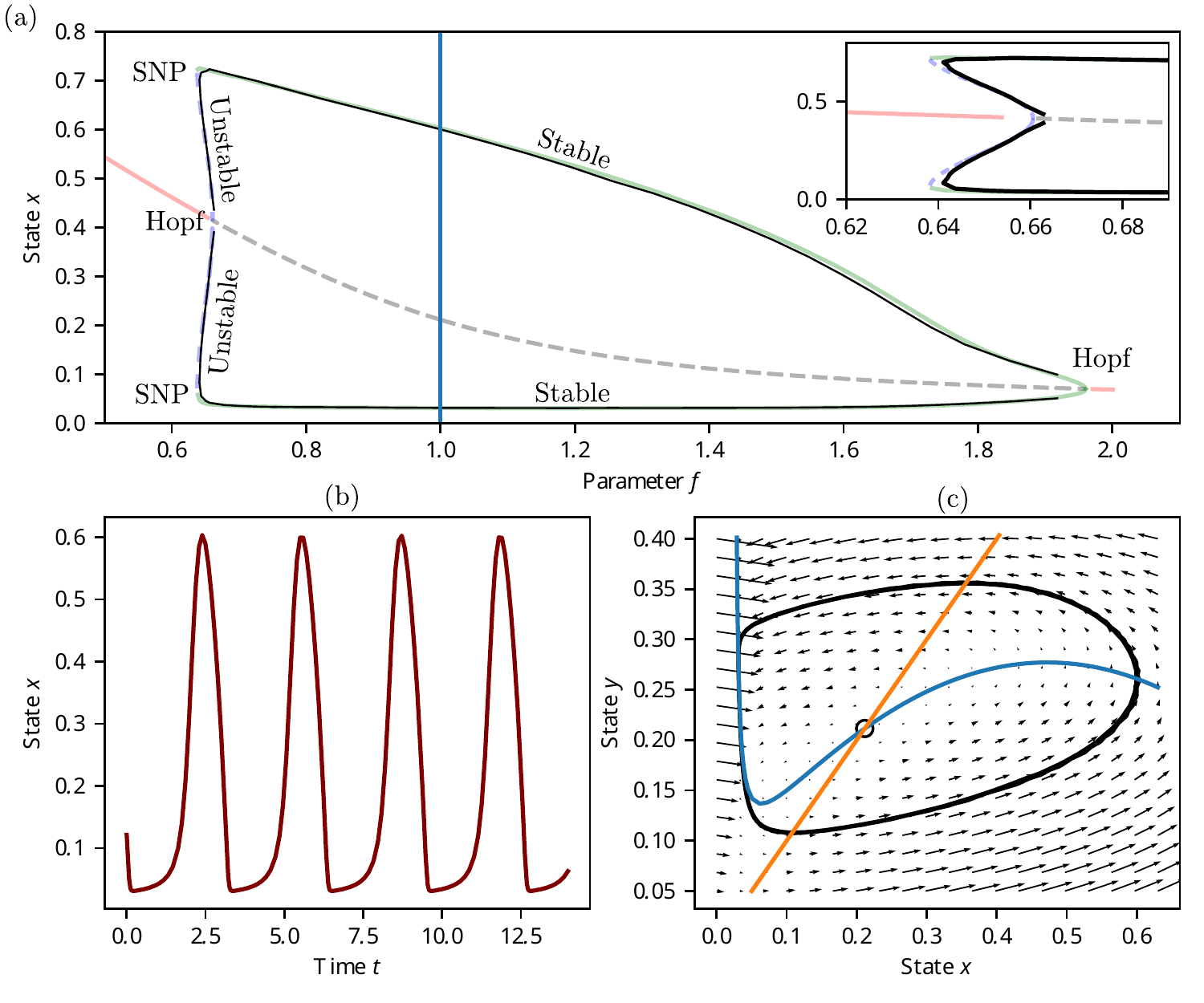}
\caption{\label{fig:org06f2023}Panel (a): bifurcation diagram of an autocatalytic chemical reaction model, described by a planar Oregonator. The dynamics of the oscillator are shown in panels (b) and (c), as a time-series plot of the \(x\) variable (panel b), and a phase-plane diagram illustrating the nullclines (blue, orange lines), limit cycle (black line), equilibrium (open circle), and flow-field (black arrows), panel (c). Time series and phase plane plots are produced at \(f = 1.0\), indicated by the vertical blue line in the bifurcation diagram. CBC results are shown (black lines) in the top panel, on top of a bifurcation diagram for the same system generated by XPPAUTO (version 8.0) with the same colour scheme as Fig. \ref{fig:org8d4d85c}. CBC continuation is terminated when the solution amplitude falls below 0.05; results are seen to match those produced by XPPAUTO.}
\end{figure*}

Simulations are produced using custom code in Python 3.
For both systems, a proportional feedback law is chosen.
A gain \(K_p=0.1\) is chosen for the gene oscillator, and gain \(K_p=4\) for the Oregonator.
Gains are chosen through trial and error, by seeking the smallest gain that stabilises the unstable branch of orbits.
Small gains are preferred for experimental applications, to avoid amplifying measurement noise.
Control is applied additively in both cases, to the state variable \(x\) of Eqs. \eqref{eq:hasty} and \eqref{eq:oregon}.
\rev{}{For both systems, the second state variable \(y\) is used as an angle-encoding embedding variable.}
Errors are calculated between an angle-encoded target \(x^*(\phi)\) and the state variable \(x\) at angle \(\phi(x,y)\).

Cubic B-splines are used for the discretisation basis functions.
A discretisation with ten coefficients is used for the gene network, and seven coefficients for the Oregonator.
Discretisation size is found by seeking the fewest coefficients \rev{to accurately describe open-loop, initialising data.}{necessary to capture open-loop, initialising data with satisfactory accuracy.}
B-spline knots are adapted after each successful prediction-correction step.
Pseudo-arclength continuation is chosen, with a secant prediction scheme.
Corrections are made by a Newton-solver, with forward finite-differences for Jacobian estimation, and a finite-differences stepsize of \(5\times10^{-3}\) for the gene oscillator, and \(1\times 10^{-2}\) for the chemical oscillator.
\rev{A maximum of three Newton iterations are taken per correction step.}{Randomness in experiments limits the accuracy of a numerical solution. It can be beneficial to cap the number of solver iterations to avoid excessive experimental effort in noisy regions. To recreate experimentally relevant conditions, we take a maximum of three Newton iterations per correction step.}
Convergence is declared when either the normed Newton-step-size or solution residual falls below \(5\times10^{-3}\), or three iterations have been reached.
Continuation is terminated when a solution amplitude falls below 0.05.

Continuation is initialised from two open-loop oscillations.
These are observed at parameters \(\gamma_y=0.03\) and \(\gamma_y=0.0301\) for the gene oscillator, and \(f=0.75\) and \(f=0.755\) for the chemical oscillator.
Two separate experiments are run, one continuing forward, then one backward.
A fixed prediction stepsize is used for the chemical oscillator, of \(0.1\) for the forward run, and \(0.05\) for the backward run.

A simple stepsize-adaption routine is chosen for the gene oscillator, to assist with continuation around the sharp fold of periodics.
Note that this method is not necessarily the best choice for CBC experiments, as it can be computationally inefficient.
Other stepsize adaptation methods for CBC are considered in \cite{schilderExperimentalBifurcationAnalysis2015}. 
Stepsize adaptation proceeds as follows.
`Current stepsize' is defined as the distance between the current corrected solution guess, and the most recently accepted solution.
The ratio between prediction stepsize and current stepsize is calculated from the solver solution at each step. 
If they are within a scale-factor of 1.2, the solution is accepted, and the stepsize is increased.
Stepsize is capped at a maximum of 0.2 for the forward-run, and 0.1 for the backward run.
If the solution distance is not within a scale-factor of 1.2, the stepsize is reduced, and a new prediction-correction step is taken.
Stepsize is capped at a minimum of \(1\times10^{-3}\).
If the stepsize is decreased below this, a convergence failure is declared, and the continuation terminates.

A polar origin must be chosen for angle-encoding, as described in section \ref{sec:org10e6e5d}.
Careful choice of polar origin is important for successful control; section \ref{sec:org6874dfd} discusses origin placement, and defines and motivates the placement options used here.
The gene oscillator uses a max-min placement for the polar origin of the forward-continuation; a min-max origin for the backward-continuation.
The chemical oscillator uses a middle-origin for the forward continuation, and a max-max origin for the backward run.
The polar origin is updated alongside the discretisation, after each new solution has been accepted.

\subsection{Results}
\label{sec:orgef114e2}
\begin{figure*}[th!]
\centering
\includegraphics[width=.8\textwidth]{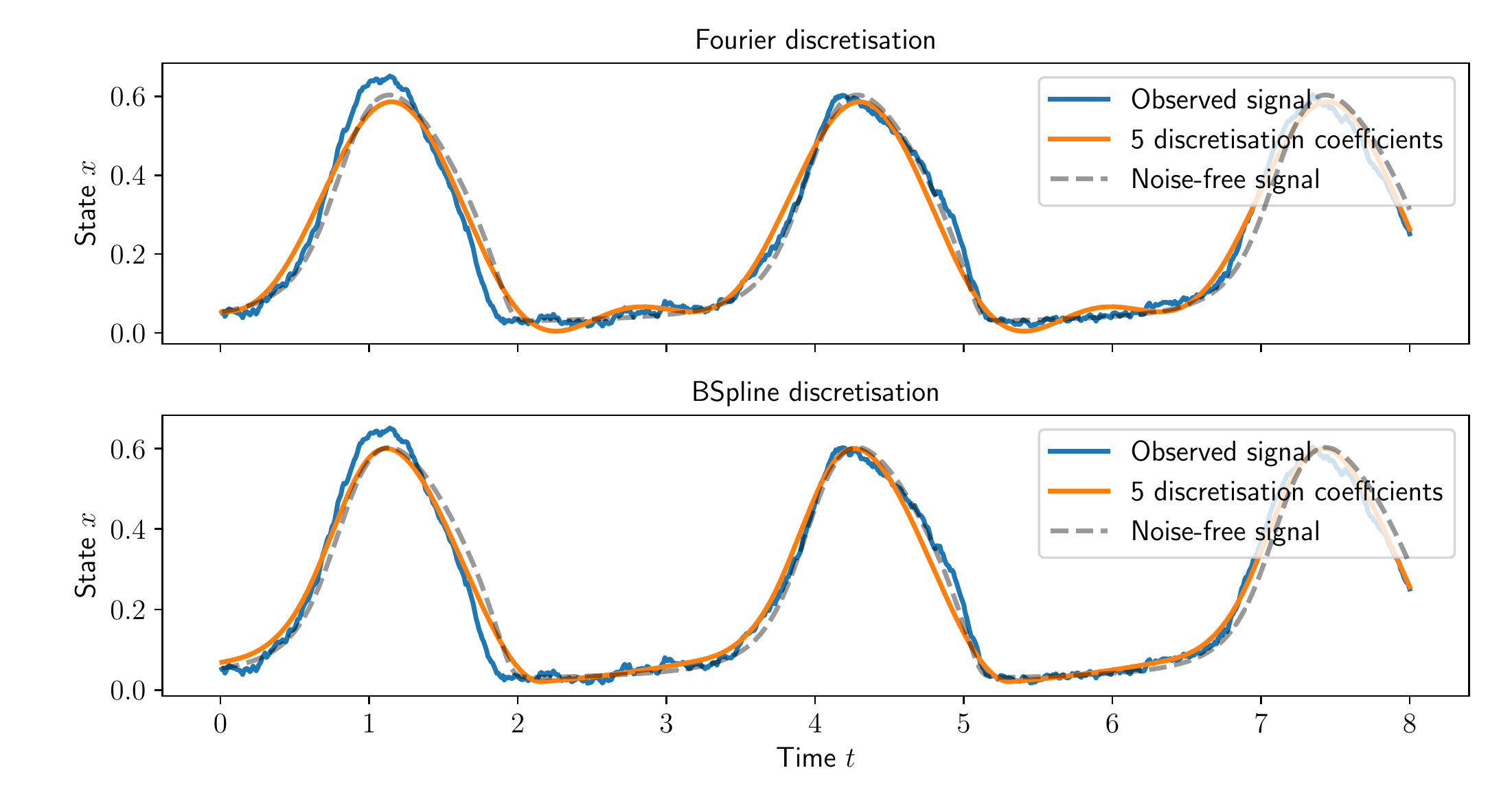}
\caption{\label{fig:org3474251}Comparison of Fourier and B-spline discretisation for noise-corrupted observations of the chemical oscillator studied in section \ref{sec:orgef114e2}, at parameter \(f=1\). A control target is generated from a deterministic, open-loop oscillation. Synthetic noisy data is generated by controlling the system using a proportional controller and the deterministically obtained control target, with controller measurements corrupted by i.i.d. Gaussian white noise of variance \(\sigma^2=0.1\), to simulate noise-corrupted observations. Spline discretisation is seen capture the true, noise-free data more accurately than Fourier.}
\end{figure*}

Gene network results are shown in Fig. \ref{fig:org8d4d85c}.
\rev{CBC results are overlaid on a bifurcation diagram generated with XPPAUTO (version 8.0) \cite{ermentroutSimulatingAnalyzingAnimating2002}.}{To validate the results, a bifurcation diagram is also generated directly from Eq.s~\eqref{eq:hasty}\eqref{eq:hasty2}, using pseudo-arclength continuation and orthogonal collocation, with XPPAUTO version 8.0.}
It is seen that CBC traces out most of the bifurcation diagram with an exact agreement to XPPAUTO.\@
\rev{}{However, unlike XPPAUTO, the CBC procedure only requires measurements of a controlled system.
This allows the system dynamics to be studied directly in an experiment, instead of through mathematical models.}
CBC is unable to trace out the canard explosion at approximately \(\gamma_y = 0.01575\).
Nevertheless, canard orbits are structurally unstable and very short-lived in planar systems \cite{dienerCanardUnchainedorHow1984}, so would be very difficult to observe in an experiment.
We therefore do not consider this a shortcoming of our methods, but rather an inherent challenge of studying physical slow-fast systems.
Note that some methods exist for controlling canard cycles \cite{durham2008feedback,jardon2021controlling}.

Fig. \ref{fig:org06f2023} shows the results of CBC on an Oregonator model.
CBC solutions are again overlaid on a bifurcation diagram generated with XPPAUTO.
Our methods are shown to trace out both stable and unstable periodic responses, with comparable accuracy to XPPAUTO.

\section{Discussion}
\label{sec:orgd7cc553}
\subsection{Robustness to observation noise}
\label{sec:org7e7811d}

Spline models are appealing for CBC: piecewise-polynomials offer the broad descriptive capabilities of polynomial fitting, without suffering from Runge's phenomenon.
In addition, the combination of model smoothness and low-degree polynomials provide compelling noise-filtering capabilities, as shown in Fig. \ref{fig:org3474251}.
Plotted are synthetic data, from a slow-fast system controlled with a proportional controller.
Data are generated from the Oregonator model described in section~\ref{sec:org6838ea7}.
To recreate a realistic experiment, controller observations are subject to measurement noise, producing a stochastic system with a noisy response.
B-spline discretisation outperforms Fourier discretisation of the noisy data.
With low discretisation sizes, Fourier discretisation struggles to match the slow-changing section of the signal; with larger discretisation sizes, it becomes susceptible to noise corruption.
B-splines provide a more accurate description of the underlying deterministic process, and are robust against observation noise.

\subsection{Angle-encoding origin choices}
\label{sec:org6874dfd}

\begin{figure*}[th!]
\centering
\includegraphics[width=\textwidth]{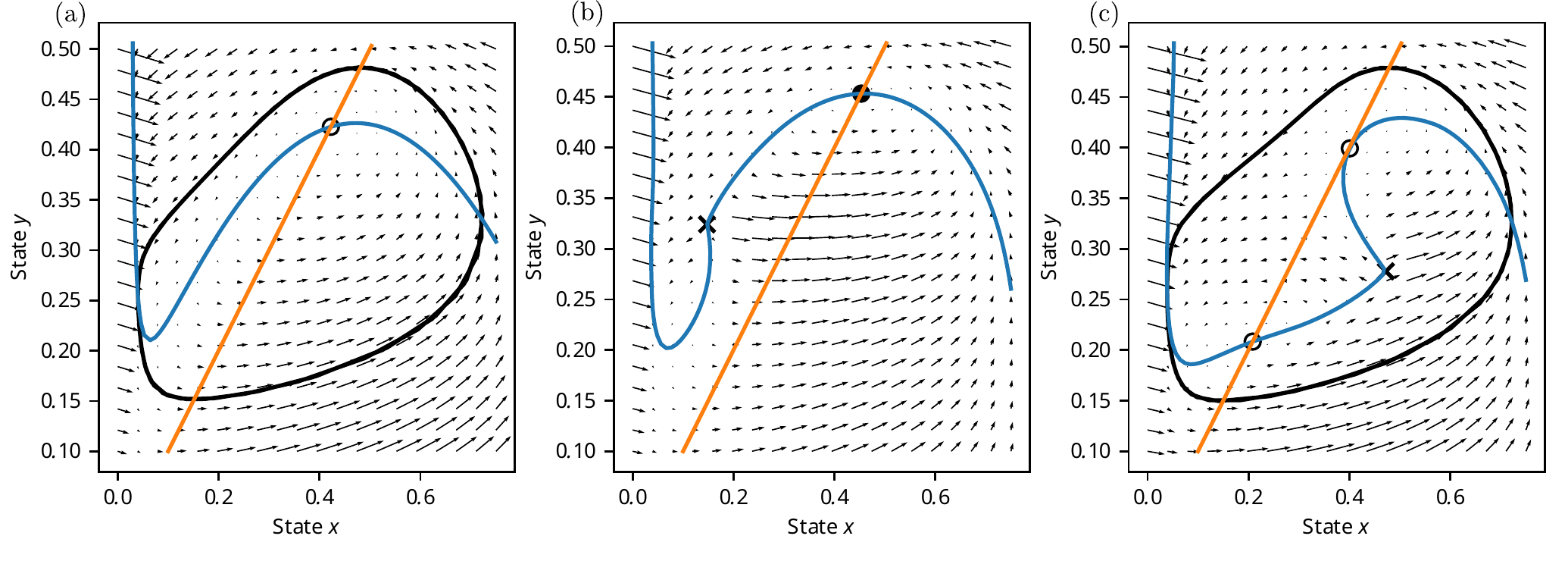}
\caption{\label{fig:org1817c06}Phase planes demonstrating mechanism of stalling for angle-encoded control, showing flow-field (arrows), \(x\) and \(y\) nullclines (blue, orange lines respectively), stable and unstable equilibria (filled, open circles), angle-encoding origin (black cross), and limit cycle trajectory where appropriate (black closed curve). Panel (a): phase-plane of an uncontrolled Oregonator at parameter \(f=0.65\). Panel (b): a noninvasive control target for parameter \(f=0.67\) is found from a discretisation of open-loop oscillations. This is then used as an angle-encoded control target at parameter \(f=0.65\), using a min-max origin. The controller is seen to put a stable equilibrium into the system, so that oscillations disappear. Panel (c): the same control target is now used with a max-min origin. Moving the polar origin has changed the nullclines, such that the previous stable equilibrium has destabilised, and the desired oscillations return.}
\end{figure*}

Angle-encoding requires a polar origin \((\mu_x, \mu_z)\), from which state angles are computed.
Easily computed choices for \(\mu_x\) and \(\mu_z\) are the mean over one period of the signals \(x(t)\) and \(z(t)\), respectively.
Nevertheless, these values do not work in all cases.
With some systems, the controller can settle to a spurious equilibrium -- an invasive equilibrium, near an equilibrium of the uncontrolled system.
Spurious equilibria are not representative of the true, uncontrolled equilibrium at that parameter value.

Fig. \ref{fig:org1817c06} shows a set of phase plane diagrams that elucidate this stalling mechanism.
Controllers modify the dynamics of a system; as angle-encoding is a time-independent control strategy, the modified dynamics can be shown on phase portrait diagrams.
The controller changes the nullcline of the controlled variable, which can introduce additional equilibria into the system, or cause changes in the position of existing equilibria.
Although noninvasive control is generally successful with angle-encoded control targets, the controlled dynamics are not always robust.
A small change in control target shape or system parameters can be enough to push one of these unstable equilibria onto the stable nullcline branch.
Such changes are to be expected within the prediction and correction steps of a continuation experiment.
Fig. \ref{fig:org1817c06} highlights this -- a noninvasive control target is found for the Oregonator, at parameter \(f=0.67\).
Periodic orbits are controlled successfully with this target, however when the same target is used at a parameter of \(f=0.65\), oscillatory dynamics are seen to have disappeared, the equilibrium has stabilised, and the system has settled to an invasive equilibrium (Fig. \ref{fig:org1817c06}, panel (b)).

In such cases the oscillatory response can often be restarted by changing \(\mu_x\) and \(\mu_z\) to move the polar origin towards the spurious equilibrium.
Changing the polar origin also changes the nullclines of the controlled system.
This can push the invasive equilibrium back onto the unstable nullcline branch, so that oscillatory dynamics reappear (Fig. \ref{fig:org1817c06}, panel (c)).
To this end, we define a selection of choices of polar origin which, while heuristic, are found to often move the nullclines in an appropriate manner as to destabilise the invasive equilibria.

\begin{figure*}[th!]
\centering
\includegraphics[width=.7\textwidth]{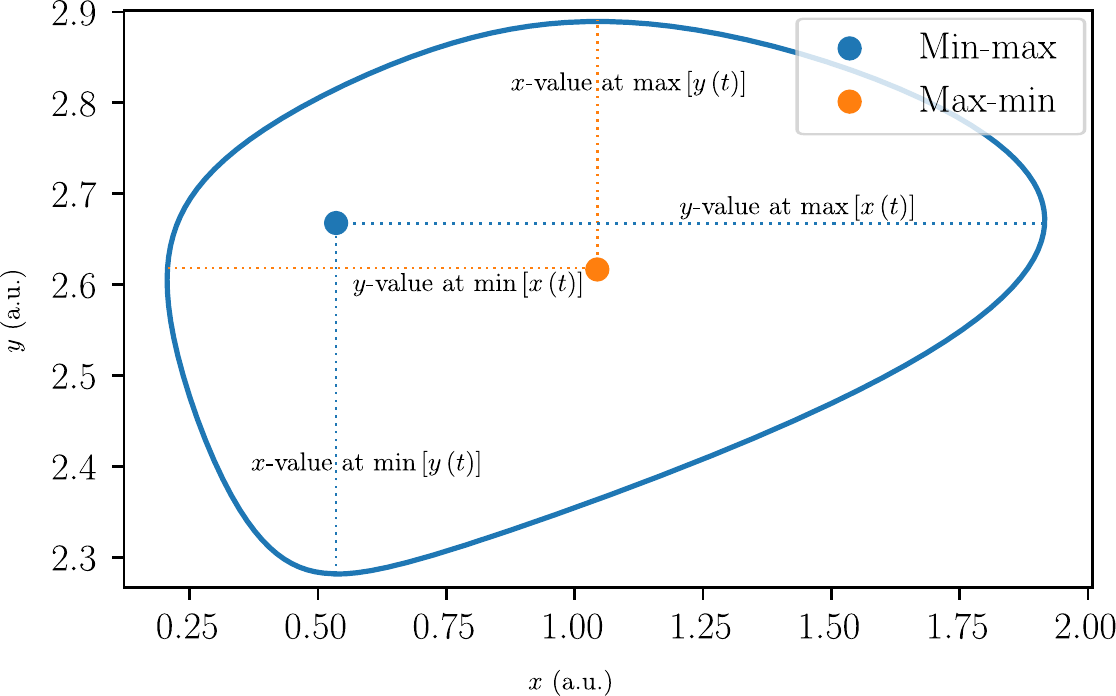}
\caption{\label{fig:org9b69e9c}Sketch of the min-max and max-min polar origin choices for determining a state angle from an embedded limit cycle.}
\end{figure*}

Two heuristics are sketched in Fig. \ref{fig:org9b69e9c}.
The first heuristic --- the min-max origin --- can fix spurious equilibria that would otherwise appear in the top-left of the reconstructed cycle, by moving the polar origin up and left slightly.
It is defined by taking \(\mu_x\) as the \(x\) coordinate of the state with minimal \(z\) value, and \(\mu_z\) as the \(z\) coordinate of the state with maximal \(x\) value.
The second heuristic --- the max-min origin --- moves the polar origin slightly down and right, and prevents spurious equilibria from appearing in the bottom-right of the reconstructed cycle.
It is defined as the \(x\) value at the state with maximal \(z\) value, and the \(z\) value of the state with minimal \(x\) value.
Additionally, a max-max and min-min origin can be defined, following the same pattern.
More alternatives are a middle-origin, where the polar origin placed in the center of a minimal bounding box for the limit cycle; to use the results of an equilibrium CBC to inform origin placement; or to manually select the polar origin and keep it fixed throughout a continuation.
These heuristics succeed in our tested cases, however it must be noted that they are merely heuristics, and may not succeed in every case, or with every embedding strategy.
Future work should investigate robust and automated polar origin placement.

\section{Conclusion}
\label{sec:org0fa60be}
CBC is a model-free method for performing a bifurcation analysis of black-box and physical systems.
To promote faster, more accurate CBC experiments for relaxation oscillations, we considered methods for adaptive discretisation using a B-spline basis.
Spline models are able to capture changeable system responses more succinctly than the standard truncated Fourier series method.
Furthermore, we proposed the angle-encoding phase constraint, which complements phase-plane CBC \cite{lee2022analysis} and produces lower-dimensional discretisations on relaxation oscillations.
Angle-encoding parameterises control targets by an angle rather than time.
It acts as a phase constraint, and removes the need to compute response frequencies.

Our angle-encoding and B-spline-discretisation methods were demonstrated on simulations of a synthetic gene oscillator \cite{hasty2002synthetic} and autocatalytic chemical oscillator \cite{tyson1980target}.
CBC produced accurate results, which visually matched those obtained using \rev{AUTO}{XPPAUTO}.\@
B-splines were found to be an efficient discretisor for CBC of slow-fast systems.
Angle-encoding offers a compelling choice of phase-constraint for autonomous systems, and both simplifies numerical computations by avoiding computation of response period, and gives robustness against drifts in response frequency.

Our methods are immediately beneficial for any system whose responses contain many Fourier harmonics.
While relaxation oscillations form the motivation of this work, \rev{the results are equally applicable to non-smooth systems, such as impact oscillators and stick-slip dynamics}{we anticipate our methods will prove beneficial to a wide range of systems. B-splines are able to succinctly capture rapidly changing responses, such as those arising in impacting and stick-slip oscillators,} and indeed any system for which Fourier discretisation is not desirable.
In addition, the phase constraint and discretisation methods are distinct.
B-spline discretisation can be used without angle-encoding on nonautonomous systems, and angle-encoding can be used with Fourier discretisation when the Fourier-Galerkin method is appropriate.
Together, these two methods open up the possibility of applying control-based continuation to a much wider range of systems --- those exhibiting autonomous responses, and those with significant amounts of higher-order harmonics.

Angle-encoding is sometimes found to cause issues within the controller.
State angles can cease to change, so that the system converges to an invasive equilibrium instead of the oscillations of interest.
Careful placement of polar origins are seen to fix this in the cases tested here.
Future work is focusing on systematically identifying when stalling happens, and methods to prevent it.

\section*{Declarations}

\subsection*{Conflict of interest}
\label{sec:org62dacf7}
The authors declare that the research was conducted in the absence of
any commercial or financial relationships that could be construed as a
potential conflict of interest.

\subsection*{Funding}
M.B. was supported by an Engineering and Physical Sciences Research Council Doctoral Training Partnership scholarship, provided by the University of Bristol. K.T.A. gratefully acknowledges the financial support of the Engineering and Physical Sciences Research Council via grant EP/T017856/1. L.M. was funded by the Engineering and Physical Sciences Research Council (EPSRC, grants EP/R041695/1 and EP/S01876X/1) and Horizon 2020 (CosyBio, grant agreement 766840). L.R. has received funding from the Royal Academy of Engineering (RF1516/15/11) which is gratefully acknowledged.

\subsection*{Acknowledgements}
The authors wish to thank Alan Champneys and Jan Sieber for helpful conversations during the preparation of this work.

\subsection*{Data availability}
Data sharing is not applicable to this article, as no datasets were generated or analysed during the current study.

\subsection*{Author Contributions}
\label{sec:orgd75e14f}
M.B. produced the research and manuscript, under the
supervision of L.M., L.R., and K.T.A.

\bibliographystyle{plain}

\begin{thebibliography}{10}

\bibitem{abeloos2022consistency}
Ga{\"e}tan Abeloos, Florian M{\"u}ller, E~Ferhatoglu, M~Scheel, Christophe
  Collette, Ga{\"e}tan Kerschen, MRW Brake, Paolo Tiso, Ludovic Renson, and
  Malte Krack.
\newblock A consistency analysis of phase-locked-loop testing and control-based
  continuation for a geometrically nonlinear frictional system.
\newblock {\em Mechanical Systems and Signal Processing}, 170:108820, 2022.

\bibitem{abeloos2023experimental}
Ga{\"e}tan Abeloos, Martin Volvert, and Ga{\"e}tan Kerschen.
\newblock Experimental characterization of superharmonic resonances using
  phase-lock loop and control-based continuation.
\newblock In {\em Nonlinear Structures \& Systems, Volume 1}, pages 131--133.
  Springer, 2023.

\bibitem{bartonControlbasedContinuationBifurcation2017}
David~AW Barton.
\newblock Control-based continuation: Bifurcation and stability analysis for
  physical experiments.
\newblock {\em Mechanical Systems and Signal Processing}, 84:54--64, 2017.

\bibitem{beregiRobustnessNonlinearParameter2021}
Sandor Beregi, David A.~W. Barton, Djamel Rezgui, and Simon~A. Neild.
\newblock Robustness of nonlinear parameter identification in presence of
  process noise using control-based continuation.
\newblock {\em arXiv:2001.11008 [nlin]}, February 2021.

\bibitem{bureauExperimentalBifurcationAnalysis2013}
Emil Bureau, Ilmar~Ferreira Santos, Jon~Juel Thomsen, Frank Schilder, and Jens
  Starke.
\newblock Experimental {{Bifurcation Analysis}} by {{Control}}-{{Based
  Continuation}}: Determining {{Stability}}.
\newblock In {\em {{ASME}} 2012 {{International Design Engineering Technical
  Conferences}} and {{Computers}} and {{Information}} in {{Engineering
  Conference}}}, pages 999--1006. {American Society of Mechanical Engineers
  Digital Collection}, September 2013.

\bibitem{chen2002model}
Luonan Chen and K.~Aihara.
\newblock A model of periodic oscillation for genetic regulatory systems.
\newblock {\em IEEE Transactions on Circuits and Systems I: Fundamental Theory
  and Applications}, 49(10):1429--1436, October 2002.

\bibitem{chia2022optogenetic}
Natalie Chia, Sang~Yup Lee, and Yaojun Tong.
\newblock Optogenetic tools for microbial synthetic biology.
\newblock {\em Biotechnology Advances}, page 107953, 2022.

\bibitem{coxNumericalEvaluationBsplines1972}
Maurice~G Cox.
\newblock The numerical evaluation of {{B}}-splines.
\newblock {\em IMA Journal of Applied Mathematics}, 10(2):134--149, 1972.

\bibitem{deboorPracticalGuideSplines1980}
Carl {de Boor}.
\newblock A {{Practical Guide}} to {{Splines}}.
\newblock {\em Mathematics of Computation}, 34(149):325, January 1980.

\bibitem{de2022control}
Irene de~Cesare, Davide Salzano, Mario di~Bernardo, Ludovic Renson, and Lucia
  Marucci.
\newblock Control-based continuation: a new approach to prototype synthetic
  gene networks.
\newblock {\em ACS Synthetic Biology}, 11(7):2300--2313, 2022.

\bibitem{denis2018identification}
Vivien Denis, M~Jossic, Christophe Giraud-Audine, B~Chomette, A~Renault, and
  Olivier Thomas.
\newblock Identification of nonlinear modes using phase-locked-loop
  experimental continuation and normal form.
\newblock {\em Mechanical Systems and Signal Processing}, 106:430--452, 2018.

\bibitem{dibernardoPredictingSyntheticGene2012}
Diego {di Bernardo}, Lucia Marucci, Filippo Menolascina, and Velia Siciliano.
\newblock Predicting synthetic gene networks.
\newblock In {\em Synthetic {{Gene Networks}}}, pages 57--81. {Springer}, 2012.

\bibitem{dienerCanardUnchainedorHow1984}
Marc Diener.
\newblock The canard unchainedor how fast/slow dynamical systems bifurcate.
\newblock {\em The Mathematical Intelligencer}, 6(3):38--49, September 1984.

\bibitem{dimatteoBayesianCurvefittingFreeknot2001}
Ilaria DiMatteo, Christopher~R Genovese, and Robert~E Kass.
\newblock Bayesian curve-fitting with free-knot splines.
\newblock {\em Biometrika}, 88(4):1055--1071, 2001.

\bibitem{durham2008feedback}
Joseph Durham and Jeff Moehlis.
\newblock Feedback control of canards.
\newblock {\em Chaos: An Interdisciplinary Journal of Nonlinear Science},
  18(1):015110, March 2008.

\bibitem{ermentroutSimulatingAnalyzingAnimating2002}
Bard Ermentrout.
\newblock {\em Simulating, Analyzing, and Animating Dynamical Systems: A Guide
  to {{XPPAUT}} for Researchers and Students}, volume~14.
\newblock {Siam}, 2002.

\bibitem{fatoorehchi2015feedback}
Hooman Fatoorehchi, Hossein Abolghasemi, Reza Zarghami, and Randolph Rach.
\newblock Feedback control strategies for a cerium-catalyzed
  belousov--zhabotinsky chemical reaction system.
\newblock {\em The Canadian Journal of Chemical Engineering}, 93(7):1212--1221,
  2015.

\bibitem{fatoorehchi2015chaos}
Hooman Fatoorehchi, Reza Zarghami, Hossein Abolghasemi, and Randolph Rach.
\newblock Chaos control in the cerium-catalyzed belousov--zhabotinsky reaction
  using recurrence quantification analysis measures.
\newblock {\em Chaos, Solitons \& Fractals}, 76:121--129, 2015.

\bibitem{field1974oscillations}
Richard~J. Field and Richard~M. Noyes.
\newblock Oscillations in chemical systems. {{IV}}. {{Limit}} cycle behavior in
  a model of a real chemical reaction.
\newblock {\em The Journal of Chemical Physics}, 60(5):1877--1884, March 1974.

\bibitem{goldschmidt1998transition}
DM~Goldschmidt, VS~Zykov, and SC~M{\"u}ller.
\newblock Transition to irregular dynamics of spiral waves under two-channel
  feedback.
\newblock {\em Physical review letters}, 80(23):5220, 1998.

\bibitem{grill1995feedback}
S~Grill, VS~Zykov, and SC~M{\"u}ller.
\newblock Feedback-controlled dynamics of meandering spiral waves.
\newblock {\em Physical review letters}, 75(18):3368, 1995.

\bibitem{grill1996spiral}
S~Grill, VS~Zykov, and Stefan~C M{\"u}ller.
\newblock Spiral wave dynamics under pulsatory modulation of excitability.
\newblock {\em The Journal of Physical Chemistry}, 100(49):19082--19088, 1996.

\bibitem{hasty2002synthetic}
Jeff Hasty, Milos Dolnik, Vivi Rottsch{\"a}fer, and James~J Collins.
\newblock Synthetic gene network for entraining and amplifying cellular
  oscillations.
\newblock {\em Physical Review Letters}, 88(14):148101, 2002.

\bibitem{jardon2021controlling}
Hildeberto Jard{\'o}n-Kojakhmetov and Christian Kuehn.
\newblock Controlling canard cycles.
\newblock {\em Journal of Dynamical and Control Systems}, pages 1--28, 2021.

\bibitem{juppApproximationDataSplines1978}
David L.~B. Jupp.
\newblock Approximation to {{Data}} by {{Splines}} with {{Free Knots}}.
\newblock {\em SIAM Journal on Numerical Analysis}, 15(2):328--343, 1978.

\bibitem{kangKnotCalculationSpline2015}
Hongmei Kang, Falai Chen, Yusheng Li, Jiansong Deng, and Zhouwang Yang.
\newblock Knot calculation for spline fitting via sparse optimization.
\newblock {\em Computer-Aided Design}, 58:179--188, January 2015.

\bibitem{khazim2021microfluidic}
Mahmoud Khazim, Elisa Pedone, Lorena Postiglione, Diego~di Bernardo, and Lucia
  Marucci.
\newblock A microfluidic/microscopy-based platform for on-chip controlled gene
  expression in mammalian cells.
\newblock In {\em Synthetic Gene Circuits}, pages 205--219. Springer, 2021.

\bibitem{kheowan2001spiral}
On-Uma Kheowan, Chi-Keung Chan, Vladimir~S Zykov, Orapin Rangsiman, and
  Stefan~C M{\"u}ller.
\newblock Spiral wave dynamics under feedback derived from a confined circular
  domain.
\newblock {\em Physical Review E}, 64(3):035201, 2001.

\bibitem{kleymanApplicationControlBasedContinuationCharacterization2020}
Gleb Kleyman, Martin Paehr, and Sebastian Tatzko.
\newblock Application of {{Control}}-{{Based}}-{{Continuation}} for
  characterization of dynamic systems with stiffness and friction
  nonlinearities.
\newblock {\em Mechanics Research Communications}, 106:103520, June 2020.

\bibitem{kleymanExperimentalApplicationControlBasedContinuation2021}
Gleb Kleyman, Martin Paehr, and Sebastian Tatzko.
\newblock Experimental {{Application}} of
  {{Control}}-{{Based}}-{{Continuation}} for {{Characterization}} of {{Isolated
  Modes}} on {{Single}}- and {{Multiple}}-{{Degree}}-of-{{Freedom Systems}}.
\newblock In {\em Nonlinear {{Structures}} \& {{Systems}}, {{Volume}} 1},
  Conference {{Proceedings}} of the {{Society}} for {{Experimental Mechanics
  Series}}, pages 135--138. {Springer International Publishing}, 2021.

\bibitem{kokotovic1999singular}
Petar Kokotovic, Hassan~K. Khali, and John O'Reilly.
\newblock {\em Singular {{Perturbation Methods}} in {{Control}}: Analysis and
  {{Design}}}.
\newblock {SIAM}, January 1999.

\bibitem{kuehn2015multiple}
Christian Kuehn.
\newblock {\em Multiple {{Time Scale Dynamics}}}.
\newblock {Springer}, February 2015.

\bibitem{kuze2019chemical}
Masakazu Kuze, Mari Horisaka, Nobuhiko~J Suematsu, Takashi Amemiya, Oliver
  Steinbock, and Satoshi Nakata.
\newblock Chemical wave propagation in the belousov--zhabotinsky reaction
  controlled by electrical potential.
\newblock {\em The Journal of Physical Chemistry A}, 123(23):4853--4857, 2019.

\bibitem{kuznetsovElementsAppliedBifurcation2013}
Yuri~A Kuznetsov.
\newblock {\em Elements of Applied Bifurcation Theory}, volume 112.
\newblock {Springer Science \& Business Media}, 2013.

\bibitem{lasalle1959time}
JP~LaSalle.
\newblock Time optimal control systems.
\newblock {\em Proceedings of the National Academy of Sciences},
  45(4):573--577, 1959.

\bibitem{lee2020reducedorder}
K.~H. Lee, D.~A.~W. Barton, and L.~Renson.
\newblock Reduced-order modelling of flutter oscillations using normal forms
  and scientific machine learning.
\newblock {\em arXiv:2011.02041 [physics]}, November 2020.

\bibitem{lee2022modelling}
KH~Lee, DAW Barton, and L~Renson.
\newblock Modelling of physical systems with a hopf bifurcation using
  mechanistic models and machine learning.
\newblock {\em arXiv preprint arXiv:2209.06910}, 2022.

\bibitem{lee2022analysis}
KH~Lee, I~Tartaruga, D~Rezgui, SA~Renson, L~Neild, and DAW Barton.
\newblock Analysis of self-excited flutter oscillations with control-based
  continuation.
\newblock {\em arXiv preprint}, 2022.

\bibitem{liHeuristicKnotPlacement2004}
Weishi Li, Shuhong Xu, Gang Zhao, and Li~Ping Goh.
\newblock A {{Heuristic Knot Placement Algorithm}} for {{B}}-{{Spline Curve
  Approximation}}.
\newblock {\em Computer-aided design and applications}, 1(1-4):727--732, 2004.

\bibitem{li2020adaptivea}
Yang Li and Harry Dankowicz.
\newblock Adaptive control designs for control-based continuation in a class of
  uncertain discrete-time dynamical systems.
\newblock {\em Journal of Vibration and Control}, 26(21-22):2092--2109,
  November 2020.

\bibitem{li2021adaptive}
Yang Li and Harry Dankowicz.
\newblock Adaptive control designs for control-based continuation of periodic
  orbits in a class of uncertain linear systems.
\newblock {\em Nonlinear Dynamics}, 103(3):2563--2579, February 2021.

\bibitem{li2022model}
Yang Li and Harry Dankowicz.
\newblock Model-free continuation of periodic orbits in certain nonlinear
  systems using continuous-time adaptive control.
\newblock {\em arXiv preprint arXiv:2203.10306}, 2022.

\bibitem{lindstromBayesianEstimationFreeknot2002}
Mary~J Lindstrom.
\newblock Bayesian estimation of free-knot splines using reversible jumps.
\newblock {\em Computational statistics \& data analysis}, 41(2):255--269,
  2002.

\bibitem{liu2003relaxation}
Weishi Liu, Dongmei Xiao, and Yingfei Yi.
\newblock Relaxation oscillations in a class of predator\textendash prey
  systems.
\newblock {\em Journal of Differential Equations}, 188(1):306--331, February
  2003.

\bibitem{lugagne2017balancing}
Jean-Baptiste Lugagne, Sebasti{\'a}n Sosa~Carrillo, Melanie Kirch, Agnes
  K{\"o}hler, Gregory Batt, and Pascal Hersen.
\newblock Balancing a genetic toggle switch by real-time feedback control and
  periodic forcing.
\newblock {\em Nature communications}, 8(1):1--8, 2017.

\bibitem{mamicAutomaticBayesianKnot2001}
G.~Mamic and M.~Bennamoun.
\newblock Automatic {{Bayesian}} knot placement for spline fitting.
\newblock In {\em Proceedings 2001 {{International Conference}} on {{Image
  Processing}} ({{Cat}}. {{No}}.{{01CH37205}})}, volume~1, pages 169--172
  vol.1, October 2001.

\bibitem{marino2011mixedmode}
F.~Marino, M.~Ciszak, S.~F. Abdalah, K.~{Al-Naimee}, R.~Meucci, and F.~T.
  Arecchi.
\newblock Mixed-mode oscillations via canard explosions in light-emitting
  diodes with optoelectronic feedback.
\newblock {\em Physical Review E}, 84(4):047201, October 2011.

\bibitem{marucciNanogDynamicsMouse2017}
Lucia Marucci.
\newblock Nanog dynamics in mouse embryonic stem cells: Results from systems
  biology approaches.
\newblock {\em Stem cells international}, 2017, 2017.

\bibitem{marucci2011derivation}
Lucia Marucci, Stefania Santini, Mario Di~Bernardo, and Diego Di~Bernardo.
\newblock Derivation, identification and validation of a computational model of
  a novel synthetic regulatory network in yeast.
\newblock {\em Journal of mathematical biology}, 62(5):685--706, 2011.

\bibitem{meijerNumericalBifurcationAnalysis2009}
Hil Meijer, Fabio Dercole, and Bart Oldeman.
\newblock Numerical {{Bifurcation Analysis}}.
\newblock In {\em Encyclopedia of Complexity and Systems Science}, pages
  6329--6352. {Springer, New York}, 2009.

\bibitem{menolascina2014vivo}
Filippo Menolascina, Gianfranco Fiore, Emanuele Orabona, Luca De~Stefano, Mike
  Ferry, Jeff Hasty, Mario Di~Bernardo, and Diego Di~Bernardo.
\newblock In-vivo real-time control of protein expression from endogenous and
  synthetic gene networks.
\newblock {\em PLoS computational biology}, 10(5):e1003625, 2014.

\bibitem{michelNewDeterministicHeuristic2021}
D.~Michel and A.~Zidna.
\newblock A new deterministic heuristic knots placement for {{B}}-{{Spline}}
  approximation.
\newblock {\em Mathematics and Computers in Simulation}, 186:91--102, August
  2021.

\bibitem{milik1996slowfast}
Alexandra Milik, Alexia Prskawetz, Gustav Feichtinger, and Warren~C. Sanderson.
\newblock Slow-fast dynamics in {{Wonderland}}.
\newblock {\em Environmental Modeling \& Assessment}, 1(1):3--17, March 1996.

\bibitem{mojrzisch2016phase}
Sebastian Mojrzisch and Jens Twiefel.
\newblock Phase-controlled frequency response measurement of a piezoelectric
  ring at high vibration amplitude.
\newblock {\em Archive of Applied Mechanics}, 86(10):1763--1769, 2016.

\bibitem{mojrzisch2012experimental}
Sebastian Mojrzisch, J{\"o}rg Wallaschek, and Jan Bremer.
\newblock An experimental method for the phase controlled frequency response
  measurement of nonlinear vibration systems.
\newblock {\em PAMM}, 12(1):253--254, 2012.

\bibitem{muller2021comparison}
Florian M{\"u}ller, Ga{\"e}tan Abeloos, Erhan Ferhatoglu, Maren Scheel,
  Matthew~RW Brake, Paolo Tiso, Ludovic Renson, and Malte Krack.
\newblock Comparison between control-based continuation and phase-locked loop
  methods for the identification of backbone curves and nonlinear frequency
  responses.
\newblock In {\em Nonlinear Structures \& Systems, Volume 1}, pages 75--78.
  Springer, 2021.

\bibitem{multamakioptogenetic}
Elina Multam{\"a}ki, Andr{\'e}s Garc{\'\i}a~de Fuentes, Oleksii Sieryi,
  Alexander Bykov, Uwe Gerken, Am{\'e}rico~Tavares Ranzani, J{\"u}rgen
  K{\"o}hler, Igor Meglinski, Andreas M{\"o}glich, and Heikki Takala.
\newblock Optogenetic control of bacterial expression by red light.

\bibitem{panagiotopoulosControlbasedContinuationUnstable2017}
Ilias Panagiotopoulos and Jens Starke.
\newblock Control-based continuation of unstable pedestrian flows.
\newblock In G\'abor St\'ep\'an and G\'abor Csern\'ak, editors, {\em
  Proceedings of 9th European Nonlinear Dynamics Conference (ENOC2017)},
  Budapest, 2017.

\bibitem{panagiotopoulos2022continuation}
Ilias Panagiotopoulos, Jens Starke, Jan Sieber, and Wolfram Just.
\newblock Continuation with non-invasive control schemes: Revealing unstable
  states in a pedestrian evacuation scenario.
\newblock {\em arXiv preprint arXiv:2203.02484}, 2022.

\bibitem{pedone2021cheetah}
Elisa Pedone, Irene De~Cesare, Criseida~G Zamora-Chimal, David Haener, Lorena
  Postiglione, Antonella La~Regina, Barbara Shannon, Nigel~J Savery, Claire~S
  Grierson, Mario Di~Bernardo, et~al.
\newblock Cheetah: a computational toolkit for cybergenetic control.
\newblock {\em ACS Synthetic Biology}, 10(5):979--989, 2021.

\bibitem{pedone2019tunable}
Elisa Pedone, Lorena Postiglione, Francesco Aulicino, Dan~L Rocca, Sandra
  Montes-Olivas, Mahmoud Khazim, Diego di~Bernardo, Maria Pia~Cosma, and Lucia
  Marucci.
\newblock A tunable dual-input system for on-demand dynamic gene expression
  regulation.
\newblock {\em Nature communications}, 10(1):1--13, 2019.

\bibitem{peralta2006controlled}
Catalina Peralta, Claudia Frank, Alex Zaharakis, Carolyn Cammalleri, Matthew
  Testa, Stephen Chaterpaul, Christian Hilaire, Daniel Lang, Daniel
  Ravinovitch, Sabrina~G Sobel, et~al.
\newblock Controlled excitations of the belousov- zhabotinsky reaction:
  experimental procedures.
\newblock {\em The Journal of Physical Chemistry A}, 110(44):12145--12149,
  2006.

\bibitem{peter2017excitation}
Simon Peter and Remco~I Leine.
\newblock Excitation power quantities in phase resonance testing of nonlinear
  systems with phase-locked-loop excitation.
\newblock {\em Mechanical Systems and Signal Processing}, 96:139--158, 2017.

\bibitem{peter2016tracking}
Simon Peter, Robin Riethm{\"u}ller, and Remco~I Leine.
\newblock Tracking of backbone curves of nonlinear systems using
  phase-locked-loops.
\newblock In {\em Nonlinear Dynamics, Volume 1}, pages 107--120. Springer,
  2016.

\bibitem{petrov1993controlling}
Valery Petrov, Vilmos Gaspar, Jonathan Masere, and Kenneth Showalter.
\newblock Controlling chaos in the belousov—zhabotinsky reaction.
\newblock {\em Nature}, 361(6409):240--243, 1993.

\bibitem{postiglione2018regulation}
Lorena Postiglione, Sara Napolitano, Elisa Pedone, Daniel~L Rocca, Francesco
  Aulicino, Marco Santorelli, Barbara Tumaini, Lucia Marucci, and Diego
  di~Bernardo.
\newblock Regulation of gene expression and signaling pathway activity in
  mammalian cells by automated microfluidics feedback control.
\newblock {\em ACS synthetic biology}, 7(11):2558--2565, 2018.

\bibitem{ren2017population}
Xinying Ren, Ania-Ariadna Baetica, Anandh Swaminathan, and Richard~M Murray.
\newblock Population regulation in microbial consortia using dual feedback
  control.
\newblock In {\em 2017 IEEE 56th annual conference on decision and control
  (CDC)}, pages 5341--5347. IEEE, 2017.

\bibitem{rensonApplicationControlbasedContinuation2019}
L~Renson, AD~Shaw, DAW Barton, and SA~Neild.
\newblock Application of control-based continuation to a nonlinear structure
  with harmonically coupled modes.
\newblock {\em Mechanical Systems and Signal Processing}, 120:449--464, 2019.

\bibitem{rensonExperimentalTrackingLimitpoint2017}
Ludovic Renson, David~AW Barton, and Simon~A Neild.
\newblock Experimental tracking of limit-point bifurcations and backbone curves
  using control-based continuation.
\newblock {\em International Journal of Bifurcation and Chaos}, 27(01):1730002,
  2017.

\bibitem{renson2016robust}
Ludovic Renson, Alicia Gonzalez-Buelga, David~AW Barton, and Simon~A Neild.
\newblock Robust identification of backbone curves using control-based
  continuation.
\newblock {\em Journal of Sound and Vibration}, 367:145--158, 2016.

\bibitem{schilderExperimentalBifurcationAnalysis2015}
Frank Schilder, Emil Bureau, Ilmar~Ferreira Santos, Jon~Juel Thomsen, and Jens
  Starke.
\newblock Experimental bifurcation analysis\textemdash{{Continuation}} for
  noise-contaminated zero problems.
\newblock {\em Journal of Sound and Vibration}, 358:251--266, 2015.

\bibitem{schwetlickLeastSquaresApproximation1995}
Hubert Schwetlick and Torsten Schutze.
\newblock Least squares approximation by splines with free knots.
\newblock {\em BIT Numerical mathematics}, 35(3):361--384, 1995.

\bibitem{seydel1991tutorial}
R.~Seydel.
\newblock Tutorial on continuation.
\newblock {\em International Journal of Bifurcation and Chaos}, 1(01):3--11,
  1991.

\bibitem{shannon2020vivo}
Barbara Shannon, Criseida~G Zamora-Chimal, Lorena Postiglione, Davide Salzano,
  Claire~S Grierson, Lucia Marucci, Nigel~J Savery, and Mario Di~Bernardo.
\newblock In vivo feedback control of an antithetic molecular-titration motif
  in escherichia coli using microfluidics.
\newblock {\em ACS Synthetic Biology}, 9(10):2617--2624, 2020.

\bibitem{showalter1979detailed}
Kenneth Showalter, Richard~M Noyes, and Harold Turner.
\newblock Detailed studies of trigger wave initiation and detection.
\newblock {\em Journal of the American Chemical Society}, 101(25):7463--7469,
  1979.

\bibitem{sieberControlBasedBifurcation2008}
Jan Sieber and Bernd Krauskopf.
\newblock Control based bifurcation analysis for experiments.
\newblock {\em Nonlinear Dynamics}, 51(3):365--377, 2008.

\bibitem{song2022bayesian}
Mingming Song, Ludovic Renson, Babak Moaveni, and Gaetan Kerschen.
\newblock Bayesian model updating and class selection of a wing-engine
  structure with nonlinear connections using nonlinear normal modes.
\newblock {\em Mechanical Systems and Signal Processing}, 165:108337, February
  2022.

\bibitem{tabor2011multichromatic}
Jeffrey~J Tabor, Anselm Levskaya, and Christopher~A Voigt.
\newblock Multichromatic control of gene expression in escherichia coli.
\newblock {\em Journal of molecular biology}, 405(2):315--324, 2011.

\bibitem{thothadri2005nonlinear}
M~Thothadri and FC~Moon.
\newblock Nonlinear system identification of systems with periodic limit-cycle
  response.
\newblock {\em Nonlinear Dynamics}, 39(1):63--77, 2005.

\bibitem{tung2002dynamics}
Chih-kuan Tung and CK~Chan.
\newblock Dynamics of spiral waves under phase feedback control in a
  belousov-zhabotinsky reaction.
\newblock {\em Physical review letters}, 89(24):248302, 2002.

\bibitem{tyson1980target}
John~J Tyson and Paul~C Fife.
\newblock Target patterns in a realistic model of the belousov--zhabotinskii
  reaction.
\newblock {\em The Journal of Chemical Physics}, 73(5):2224--2237, 1980.

\bibitem{valenzuelaEvolutionaryComputationOptimal2013}
Olga Valenzuela, Blanca {Delgado-Marquez}, and Miguel Pasadas.
\newblock Evolutionary computation for optimal knots allocation in smoothing
  splines.
\newblock {\em Applied Mathematical Modelling}, 37(8):5851--5863, 2013.

\bibitem{vernon1968relaxation}
F.~L. Vernon and R.~J. Pedersen.
\newblock Relaxation {{Oscillations}} in {{Josephson Junctions}}.
\newblock {\em Journal of Applied Physics}, 39(6):2661--2664, May 1968.

\bibitem{virtanenSciPyFundamentalAlgorithms2020}
Pauli Virtanen, Ralf Gommers, Travis~E. Oliphant, Matt Haberland, Tyler Reddy,
  David Cournapeau, Evgeni Burovski, Pearu Peterson, Warren Weckesser, Jonathan
  Bright, St{\'e}fan~J. {van der Walt}, Matthew Brett, Joshua Wilson, K.~Jarrod
  Millman, Nikolay Mayorov, Andrew R.~J. Nelson, Eric Jones, Robert Kern, Eric
  Larson, C.~J. Carey, {\.I}lhan Polat, Yu~Feng, Eric~W. Moore, Jake
  VanderPlas, Denis Laxalde, Josef Perktold, Robert Cimrman, Ian Henriksen,
  E.~A. Quintero, Charles~R. Harris, Anne~M. Archibald, Ant{\^o}nio~H. Ribeiro,
  Fabian Pedregosa, and Paul {van Mulbregt}.
\newblock {{SciPy}} 1.0: Fundamental algorithms for scientific computing in
  {{Python}}.
\newblock {\em Nature Methods}, 17(3):261--272, March 2020.

\end{thebibliography}

\end{document}